\documentclass[11pt]{amsart}

\usepackage[margin=1.08in]{geometry}
\usepackage{amsmath,amssymb,amsthm,mathtools,mathrsfs}
\usepackage{booktabs,array}
\usepackage[expansion=false]{microtype}
\usepackage{enumitem}
\usepackage[hidelinks]{hyperref}
\usepackage{xcolor}

\newtheorem{theorem}{Theorem}[section]
\newtheorem{proposition}[theorem]{Proposition}

\newtheorem{corollary}[theorem]{Corollary}
\theoremstyle{remark}
\newtheorem{remark}[theorem]{Remark}

\theoremstyle{plain}
\newtheorem*{cwfourfold}{Theorem (integral Cohen--Wolfart--TTV comparison)}
\newtheorem*{cartanlevels}{Proposition (Cartan-normalizer levels)}
\theoremstyle{remark}
\newtheorem*{newremark}{Remark}

\newcommand{\Q}{\mathbf Q}
\newcommand{\F}{\mathbf F}
\newcommand{\Z}{\mathbf Z}
\newcommand{\C}{\mathbf C}
\newcommand{\HH}{\mathfrak H}
\newcommand{\OK}{\mathcal O_K}
\newcommand{\Pone}{\mathbf P^1}
\newcommand{\PSL}{\mathrm{PSL}_2}
\newcommand{\SL}{\mathrm{SL}_2}
\newcommand{\Gal}{\mathrm{Gal}}
\newcommand{\tr}{\operatorname{tr}}
\newcommand{\Norm}{\operatorname{Norm}}

\newcommand{\Aut}{\operatorname{Aut}}

\title[A hypergeometric RM surface and level-11 covers]
{A Hypergeometric RM Surface: Integral Monodromy and Level-11 Covers}
\author{Devendra Tiwari}
\date{November 2026}

\begin{document}

\begin{abstract}
We study the Jacobians of a genus-two family with real multiplication by
$K=\Q(\sqrt5)$. We realize these abelian surfaces as factors of
cyclic-cover Jacobians and determine their integral rank-two monodromy.
For the real embedding with $\sqrt5>0$, the projective image is the Hecke
group $\Delta(2,5,\infty)$; the conjugate real image is nondiscrete.
After a finite base change, we give an explicit degree-four map from
Cohen--Wolfart's Hecke cyclic curve to the genus-two curve. Four deck translates of the induced
Jacobian pullback give an integral isomorphism from the fourth power of
the surface to the eight-dimensional primitive order-$20$ factor, with
an explicit nonproduct polarization. At the primes above $11$, the
degree-$12$ torsion-line covers are normalizations of quadratic fibre
products of genus-one triangular modular curves. We also describe the
associated split and nonsplit Cartan-normalizer covers. Finite-field
hypergeometric traces give the local Euler factors and a quadratic
transformation with a twist by the quadratic character of $2$.
\end{abstract}

\maketitle

\section{Introduction and main results}

\subsection{From the Legendre family to hypergeometric abelian varieties}

The Legendre family gives the classical model for the passage from an
algebraic family to periods, monodromy and level structures.  Its periods
satisfy a Gauss hypergeometric equation, integral homology supplies the
lattice on which monodromy acts, and reduction of that lattice produces the
usual level covers of the parameter line.  Hypergeometric abelian
varieties extend this picture to higher-dimensional Jacobians, but with an
important extra step: the rank-two differential system may live on a proper
abelian factor of a cyclic-cover Jacobian rather than on the whole Jacobian.
Thus the factor carrying the periods and its integral homology must be
identified before reduction modulo a prime of the real-multiplication field
has an intrinsic geometric meaning.

A convenient algebraic realization starts from generalized Legendre curves
of the form
\[
 y^N=x^A(1-x)^B(1-tx)^C,
 \qquad t\in\Pone\setminus\{0,1,\infty\}.
\]
The character eigenspaces give Gauss hypergeometric local systems.
Rational Galois-stable combinations of them define abelian factors
\cite{Archinard2003,Darmon1999}. Our example is the genus-two family
of Tautz, Top and Verberkmoes (TTV), whose Jacobian has real
multiplication (RM) \cite{TTV1991}. We first identify this surface
inside the primitive factor of the Jacobian of a degree-ten cyclic
cover.

\subsection{Hypergeometric monodromy, triangle groups, and Darmon's dictionary}

To keep two different kinds of triples separate, we write
$(\alpha,\beta;\gamma)$ for the parameters of
${}_2F_1(\alpha,\beta;\gamma;t)$ and
$(e_0,e_1,e_\infty)$ for a triangle signature.  The local exponent
differences of the Gauss equation at $0,1,\infty$ are
\[
 1-\gamma,\qquad \gamma-\alpha-\beta,\qquad \beta-\alpha.
\]
A nonintegral rational exponent difference gives local projective
monodromy whose order is its denominator in lowest terms. In the systems
considered here, zero exponent differences give nontrivial parabolic
monodromy, recorded by $\infty$. These local orders determine an
abstract triangle orbifold group through which the projective monodromy
factors. They do not determine its global image or establish that a
particular real representation is faithful or discrete.

In his fourteenth lecture, Darmon formulates a geometric dictionary for
generalized Fermat equations in which the relevant rigid families are
organized by triangle-group signatures \cite[Principle~6]{Darmon2004}.  Up to permutation of
the three singular points, the four one-parameter branches lead in
characteristic zero to the types
\[
 (\infty,\infty,\infty),\qquad
 (r,\infty,\infty),\qquad
 (r,r,\infty),\qquad
 (q,r,\infty),
\]
and the associated finite covers arise from torsion in the corresponding
hypergeometric abelian varieties
\cite{Darmon2000}.  The present paper concerns the
$(r,\infty,\infty)$ branch with $r=5$ and studies the geometric variation of
the family before imposing arithmetic conditions on individual
specializations.

When the relevant real projective monodromy is a discrete triangle group,
there is a second description of its finite levels.  Congruence subgroups of
triangle groups, defined through their trace fields and quaternionic orders,
give algebraic three-point covers, called triangular modular curves
\cite{ClarkVoight2019,DuqueRoseroVoight2023}. A Borel subgroup stabilizes
a line in a two-dimensional residual representation; its quotient cover
therefore records the choice of such a line.  In this sense triangular
modular curves play for hypergeometric abelian varieties a role analogous to
modular curves with level structure for elliptic curves.  The point of the present
example is that these two constructions meet only after the correct abelian
factor and integral monodromy have been fixed.

The original genus-two construction and the Hecke-period interpretation
also admit a concrete integral comparison with Cohen--Wolfart's
\emph{Hecke} example \cite[\S6.1]{CohenWolfart1990}.  Over an explicit
finite base change we obtain a degree-four morphism from their order-$20$
cyclic curve to $C_t^-$; four deck translates of its Jacobian pullback
give an isomorphism from the fourth power of the surface to the
eight-dimensional primitive factor, with nonzero off-diagonal
polarization pairings.  This is distinct from Cohen--Wolfart's
other $(5,\infty,\infty)$ example, whose marked projective representations
are inequivalent to ours (Theorem~\ref{thm:C}).  The abstract fourth-power
\emph{isogeny} was already described by Cohen--W\"ustholz
\cite[pp.~98--99]{CohenWustholz2002}; the intrinsic primitive polarization
type and a related genus-five splitting also have antecedents
\cite{WangXueYu2014,LangeRodriguezRojas2014}.
The Hilbert-modular realization of this TTV family likewise belongs to
the prior literature \cite[\S5.2.1]{BogoLiAtkin2026}.

\subsection{The genus-two RM family and its level geometry}

We study the family
\begin{equation}\label{eq:TTVcurve}
 C_t^-:\qquad Y^2=X^5-5X^3+5X+2-4t,
 \qquad t\in U:=\Pone\setminus\{0,1,\infty\},
\end{equation}
whose Jacobian $A_t=J(C_t^-)$ has real multiplication by
\begin{equation}\label{eq:K}
 K=\Q(\sqrt5)
\end{equation}
\cite[Theorem~1]{TTV1991}. Write $\sigma_+$ for the real embedding with
$\sigma_+(\sqrt5)>0$ and $\sigma_-$ for its conjugate. The two embeddings
give the rank-two Picard--Fuchs systems
\begin{equation}\label{eq:TTVPF}
 {}_2F_1\!\left(\frac3{10},\frac7{10};1;t\right),
 \qquad
 {}_2F_1\!\left(\frac1{10},\frac9{10};1;t\right),
\end{equation}
with the RM period interpretation made explicit in
\cite[Theorems~4.2.4--4.2.5]{Takei2014}.  Their exponent differences are
$(0,0,-2/5)$ and $(0,0,-4/5)$, so both have local projective orders
$(\infty,\infty,5)$.  This is the expected local type for the $r=5$ case of
Darmon's $(r,\infty,\infty)$ branch. Under $\sigma_+$ the projective
image is $\Delta(2,5,\infty)$, whereas under $\sigma_-$ it is nondiscrete.
The Hecke image is additional information, not a consequence of the
three local orders alone.

Let $\Lambda_t=H_1(A_t,\Z)$, with its rank-two $\mathcal O_K$-module
structure. Reduction of this integral homology lattice defines the
level-$11$ construction. For a coefficient prime $\mathfrak p\mid11$,
\[
 A_t[\mathfrak p](\C)\simeq\Lambda_t/\mathfrak p\Lambda_t,
 \qquad
 \mathbf P(A_t[\mathfrak p])\simeq\Pone(\F_{11}).
\]
A point of this projective line is a one-dimensional $\F_{11}$-subspace,
equivalently a cyclic subgroup of order $11$ in the $\mathfrak p$-torsion
component.  Its stabilizer is a Borel subgroup of $\PSL(\F_{11})$, so the
associated degree-$12$ cover will be called the Borel level cover of the
family.  This is a relative level interpretation for the given RM family; no
universal fine moduli statement for all hypergeometric abelian varieties is
needed.

On the triangle-group side, the same Borel construction gives level-$11$
triangular modular curves for $\Delta(2,5,\infty)$.  The two descriptions are
related by the classical quadratic parameter
\begin{equation}\label{eq:quadratic-map}
 s=4t(1-t).
\end{equation}
The central statement of Section~\ref{sec:level11} identifies the covers
of $A_t$ with the connected normalizations of the fibre products of
these Hecke covers under \eqref{eq:quadratic-map}. This is what we mean
by a \emph{normalized quadratic pullback}. Throughout the comparison,
a \emph{marked cover} has its branch values ordered as $(0,1,\infty)$
on the fixed parameter line; its sheets are not individually labelled.
The positively oriented based loops used to compare monodromy are fixed
in Section~\ref{sec:monodromy}.  Because \eqref{eq:quadratic-map} is
ramified at the order-five point, this is a pullback of projective structures,
not the ordinary orbifold covering associated with an index-two subgroup.

The same family also appears on the arithmetic side of Darmon's programme.
Modern work relates a TTV-form model to the adjacent odd $(r,r,p)$ branch
after explicit parameter change, base extension and quadratic twisting
\cite{BCDFI2025}; the original parameter lines and the finite covers studied
here are not thereby identified.  A companion paper treats that adjacent
$(5,5,\infty)$ geometry \cite{TiwariDarmonCompanion2026}.  Finally, the
cyclic-cover realization has compatible $\ell$-adic and finite-field
realizations \cite{DFLST2016,FuselierLongRamakrishnaSwisherTu2022}.  Section
\ref{sec:finite-field} uses them to compute the Frobenius traces and local
Euler factors of the same two RM constituents and to identify the arithmetic
quadratic transformation, including its constant character twist.

At the chosen coefficient prime $\mathfrak p_1\mid11$, the same torsion local system
also yields arithmetic split and nonsplit Cartan-normalizer covers over
$K$: these record Galois-stable pairs of projective torsion lines rather
than individual cyclic subgroups.  Their degrees and genera, and a
special CM fibre, illustrate the arithmetic scope of the level geometry,
without determining their rational points or proving a uniform image
theorem (Section~\ref{sec:cartan11}).

\subsection{Main results}

The first main result makes the relation between the RM surface and Archinard's
cyclic-cover construction explicit in the isogeny category.  Consider
\begin{equation}\label{eq:X-Y-cyclic}
 \mathcal X_t:\ y_X^{10}=x^3(1-x)^7(1-tx)^3,
 \qquad
 \mathcal Y_t:\ y_Y^{10}=x(1-x)^9(1-tx).
\end{equation}
Their smooth projective models are isomorphic over $\Q(t)$.  TTV's common
hyperelliptic model is
\begin{equation}\label{eq:H}
 H_t:\qquad v^2=u^{20}+(2-4t)u^{10}+1.
\end{equation}
Write $J^{\mathrm{new}}(\mathcal X_t)$ for the identity component of
the intersection of the kernels of the maps to Jacobians of lower-degree
cyclic quotients. This primitive factor has dimension four; its construction
is recalled in Section~\ref{sec:factor}. Real multiplication
decomposes the first cohomology of $A_t$, over $\mathbf R$, into two rank-two
local systems indexed by the real embeddings of $K$. Their Hodge
filtrations make them variations of Hodge structure of weight one.

\begin{theorem}[Explicit realization of the RM surface]\label{thm:A}
There is an abelian surface $B_t$ over $\Q(t)$ and an isogeny
\begin{equation}\label{eq:factor-main}
 J^{\mathrm{new}}(\mathcal X_t)
 \sim_{\Q(t)} J(C_t^-)\times B_t.
\end{equation}
The first factor is induced by the degree-four morphism
\begin{equation}\label{eq:mmap}
 m:H_t\longrightarrow C_t^- ,\qquad
 (X,Y)=\left(u^2+u^{-2},\frac{v}{u^5}\right).
\end{equation}
After complexification, the pullbacks of the two rank-two variations
lie in the sums of character local systems indexed by $\{3,7\}$ and
$\{1,9\}$, respectively. Projection to either character summand in the
corresponding pair is an isomorphism. The abelian factor over $\Q(t)$
is selected by a rational idempotent built from the additional involution.  Over
$L=\Q(\zeta_5)$ one further has
\begin{equation}\label{eq:square}
 B_{t,L}\sim J(C_t^-)_{L},
 \qquad
 J^{\mathrm{new}}(\mathcal X_t)_L\sim J(C_t^-)_{L}^{\,2}.
\end{equation}
\end{theorem}

The common cover and the dihedral splitting mechanism have antecedents in
TTV; Theorem~\ref{thm:A} records the rational idempotent and
explicit identification of the differentials with Archinard's primitive factor.  Darmon's
hypergeometric abelian-variety realization provides the broader context
\cite[Proposition~3.2]{Darmon1999}.

Put
\begin{equation}\label{eq:alpha-mu}
 \alpha=\frac{\sqrt5-1}{2},\qquad
 \OK=\Z[\alpha],\qquad
 \mu=2+\alpha=\frac{3+\sqrt5}{2}.
\end{equation}
The following proposition gives the monodromy matrices in an integral
basis of the homology lattice.

\begin{proposition}[Integral RM monodromy]\label{prop:B}
There is an $\OK$-basis of $\Lambda_t=H_1(A_t,\Z)$ for which
\begin{equation}\label{eq:exact-monodromy}
 M_0=\begin{pmatrix}1&1\\0&1\end{pmatrix},\qquad
 M_1=\begin{pmatrix}1&0\\-\mu&1\end{pmatrix},\qquad
 M_\infty=\begin{pmatrix}1&-1\\\mu&1-\mu\end{pmatrix},
\end{equation}
with $M_0M_1M_\infty=I$ and $M_\infty^5=-I$. Under $\sigma_+$
the projective image is $\Delta(2,5,\infty)$; under $\sigma_-$
it is nondiscrete.
\end{proposition}

Proposition~\ref{prop:B} records the integral normalization that will be
used in the comparison with Cohen--Wolfart and in reduction modulo $11$.
Cohen--Wolfart's \emph{different} Example~1 on pp.~104--105 instead
uses a discrete uniformizing group $\Delta(5,\infty,\infty)$
\cite[Example~1]{CohenWolfart1990}.  Their separate Hecke construction on p.~106 admits the period-map
comparison in Section~\ref{sec:cw-hecke}.
Takei
explicitly raised the question of comparing the two constructions
\cite[Chapter~6, item~2]{Takei2014}.

\begin{theorem}[Comparison with Cohen--Wolfart's Example~1]\label{thm:C}
After fixing the conventions used in Cohen--Wolfart's Example~1, its
two projective monodromy representations, indexed by the real embeddings
of the RM field, are not conjugate to the two representations of the
TTV family under
the allowed permutations of cusps, scalar gauges, duality, exchange of real
embeddings, Galois conjugation or orientation reversal.  Consequently, the
specified TTV period map and the modular embedding of Cohen--Wolfart's
Example~1 cannot be identified by these operations.
\end{theorem}

The obstruction is already visible in the two parabolic monodromies.  In the
normal form
\[
 A=\begin{pmatrix}1&1\\0&1\end{pmatrix},\qquad
 B_u=\begin{pmatrix}1&0\\-u&1\end{pmatrix},
\]
one has
\begin{equation}\label{eq:kappa}
 \kappa:=\tr[A,B_u]-2=u^2.
\end{equation}
The $u$-values for the RM family form the Galois orbit defined by $X^2-3X+1$, while the
Cohen--Wolfart orbit is defined by $X^2-5X+5$; their squared-value orbits are
disjoint.

The level-$11$ result concerns the Borel level at the two primes above
$11$.  Once the geometric lattice in Proposition~\ref{prop:B} is fixed,
reduction gives a two-dimensional $\F_{11}$-space.  Its projective space
$\Pone(\F_{11})$ has $12$ points, and the Borel quotient records a
one-dimensional subspace of $A_t[\mathfrak p]$, equivalently a cyclic
subgroup of order $11$ in that RM torsion component.  On the Hecke side, the same
line-stabilizer construction is the standard Borel level of the triangular
congruence theory of Clark--Voight; the two low-genus possibilities at this
level occur in the tables of Duque--Rosero--Voight
\cite{ClarkVoight2019,DuqueRoseroVoight2023}.  Our task is to show that the
level obtained from $A_t$ is the normalized quadratic pullback of
that Hecke level with the peripheral generators kept in the same order.

Write
\begin{equation}\label{eq:split11}
 11\OK=\mathfrak p_1\mathfrak p_2,
 \qquad
 \mathfrak p_1=(11,\alpha-3),\quad
 \mathfrak p_2=(11,\alpha-7),
\end{equation}
and let $G=\PSL(\F_{11})$.  The projective reductions of
\eqref{eq:exact-monodromy} are surjective onto $G$.  Their degree-$12$ Borel
quotients have the common ordered passport
\begin{equation}\label{eq:passport}
 (11\,1;\ 11\,1;\ 5^2 1^2)
\end{equation}
and genus $3$, but the two ordered generating triples remain inequivalent.
The relation with the Hecke model uses the known quadratic parameter
\begin{equation}\label{eq:quadratic-map-intro}
 \psi:\Pone_t\longrightarrow\Pone_s,
 \qquad s=4t(1-t).
\end{equation}
The quadratic transformation of the Gauss equation and the Hecke
interpretation are known. The statement below identifies the finite
covers obtained by reducing the integral homology representation.

\begin{theorem}[Level-$11$ quadratic pullback]\label{thm:E}
For each $\mathfrak p_i\mid11$, let
$q_i:\Delta(2,5,\infty)\twoheadrightarrow G$ be the reduction of the
integral monodromy representation, let $Z_i\to\Pone_s$ be the corresponding
Galois (regular) Hecke cover, and put $E_i=Z_i/\mathcal B$ for a Borel subgroup
$\mathcal B<G$.  Then, geometrically over $\C$ and with the branch points ordered as above,
\begin{equation}\label{eq:regular-fibre-product-intro}
 X_i^{\mathrm{proj}}\simeq
 \Norm\!\left(\Pone_t\times_{\Pone_s}Z_i\right),
\end{equation}
and
\begin{equation}\label{eq:borel-fibre-product-intro}
 Y_i\simeq
 \Norm\!\left(\Pone_t\times_{\Pone_s}E_i\right).
\end{equation}
The regular Hecke cover has degree $660$, inertia $(11,2,5)$ and genus
$70$; the Borel quotient has passport $(11\,1;2^6;5^2 1^2)$ and genus $1$.
The degree-two maps $X_i^{\mathrm{proj}}\to Z_i$ and $Y_i\to E_i$ are
branched at exactly $132$ and $4$ points above $s=\infty$, respectively,
and recover the genera $205$ and $3$.
\end{theorem}

The scalar quadratic transformation, the Hecke group, and Darmon's finite
projective realization are known ingredients.  What Theorem~\ref{thm:E}
adds in the present setting is the identification of the reductions of the
integral monodromy of $A_t$ with the corresponding Hecke branch cycles and the
resulting connected normalized fibre products.

There is also a useful arithmetic contrast on the classical modular side.
Under $U_t\simeq\bar\Gamma(2)\backslash\HH$, the regular and Borel-preimage
subgroups have generalized level $110$ but no cusp of width $110$.  The
Fricke--Wohlfahrt level theorem together with Larcher's cusp-amplitude criterion
therefore shows that both are noncongruence subgroups of
$\bar\Gamma(2)$.  We record this as Corollary~\ref{cor:noncongruence} after
the geometric pullback theorem has been proved.

\subsection{Organization of the paper}

Section~\ref{sec:factor} establishes the algebraic realization of the RM
surface inside the primitive factor of the degree-ten generalized-Legendre
cover's Jacobian;
Theorem~\ref{thm:A} is proved in the explicit factorization
subsection, Section~\ref{sec:factor-explicit}.  Section~\ref{sec:monodromy}
then determines the geometric $\mathcal O_K$-lattice and the two real
projective images; Proposition~\ref{prop:B} is proved in
Section~\ref{sec:integral-lattice}, and the comparison in
Theorem~\ref{thm:C} is proved in Section~\ref{sec:cw-comparison}.

Section~\ref{sec:level11} reduces that lattice at the two primes above $11$,
constructs the regular and Borel level covers, and compares them with the
triangular modular levels of the Hecke group.  It also proves a concrete
integral fourfold comparison with Cohen--Wolfart's Hecke cyclic curve and
describes the arithmetic Cartan-normalizer covers.  Theorem~\ref{thm:E} culminates
in the normalized fibre-product statements of Section~\ref{sec:quadratic-pullback};
the $\Q(\sqrt5)$-models of the genus-one Borel curves and their genus-three
pullbacks are established in Section~\ref{sec:borel-descent}.  The same
section also records the logically separate classical noncongruence result.
Section~\ref{sec:finite-field} gives the arithmetic realization of the two
rank-two systems: Proposition~\ref{prop:finite-field-euler} is proved in
Section~\ref{sec:frob-euler}, and the quadratic transformation of Frobenius
traces is established in Section~\ref{sec:frob-quadratic}.  The final section
collects the remaining questions on explicit equations, arithmetic
specializations and moduli.
\section{The RM surface inside a generalized-Legendre Jacobian}
\label{sec:factor}

The level-$11$ covers are attached to $A_t=J(C_t^-)$ and its integral
homology. We first locate this surface in the primitive factor of the
Jacobian of a degree-ten generalized-Legendre cover. TTV supplies the
auxiliary cover and RM quotients; Archinard's primitive-factor
construction organizes the relevant characters. The maps and rational
idempotents below identify the abelian factor carrying the two period
systems in \eqref{eq:TTVPF}. Section~\ref{sec:monodromy} identifies
its integral homology lattice; reduction of this lattice in
Section~\ref{sec:level11} then gives the torsion representations.

\subsection{Auxiliary geometry of the RM family}

The real-multiplication family \eqref{eq:TTVcurve} is obtained in TTV from a
genus-five auxiliary curve \cite[Proposition~3 and the following remark]{TTV1991}.  Put $c=2-4t$ and
\begin{equation}\label{eq:D}
 D_t^-:\qquad w^2=z(z^{10}+cz^5+1).
\end{equation}
The automorphisms relevant here are
\begin{equation}\label{eq:r-sigma}
 r(z,w)=(\zeta_5z,\zeta_5^3w),
 \qquad
 \sigma(z,w)=(z^{-1},w/z^6),
\end{equation}
with $\sigma r\sigma=r^{-1}$.  The quotient by $\sigma$ is
$C_t^-$ via
\begin{equation}\label{eq:D-C-map}
 \phi:D_t^-\longrightarrow C_t^- ,\qquad
 (X,Y)=\left(z+z^{-1},\frac{w}{z^3}\right).
\end{equation}
These formulas are part of the explicit TTV/Takei geometry
\cite{TTV1991,Takei2014}.

For
\[
\eta_j=\frac{z^j\,dz}{w},\qquad 0\le j\le4,
\]
one has
\begin{equation}\label{eq:diff-action}
 r^*\eta_j=\zeta_5^{j-2}\eta_j,
 \qquad
 \sigma^*\eta_j=-\eta_{4-j}.
\end{equation}
In particular
\begin{equation}\label{eq:pullback-Cdiffs}
 \phi^*\left(\frac{dX}{Y}\right)=\eta_3-\eta_1,
 \qquad
 \phi^*\left(\frac{X\,dX}{Y}\right)=\eta_4-\eta_0.
\end{equation}
These are the two RM eigendifferentials whose periods satisfy
\eqref{eq:TTVPF}.

There is also an elliptic quotient
\begin{equation}\label{eq:q-ell}
 q:D_t^-\longrightarrow E_t,
 \qquad
 (s,V)=(z^5,wz^2),
\end{equation}
with
\[
 E_t:\qquad V^2=s(s^2+cs+1).
\]
Let
\begin{equation}\label{eq:Pdef}
 P=\ker(q_*)^0\subset J(D_t^-).
\end{equation}
Then $\dim P=4$ and its holomorphic differentials are spanned by
$\eta_0,\eta_1,\eta_3,\eta_4$.

\subsection{The primitive degree-ten factor}

In the notation of \cite[Sec.~6.2]{FuselierLongRamakrishnaSwisherTu2022},
our curve $\mathcal Y_t$ is the generalized Legendre curve
$C_t^{[10;1,9,1]}$.  For every proper divisor $d\mid10$ there is a natural
map to the degree-$d$ cyclic cover, and the primitive Jacobian is the identity
component of the intersection of the kernels of the induced maps on
Jacobians.  We use $J^{\mathrm{new}}$ for this primitive part, as in Archinard's
construction.  The formulation is useful here because it separates
the rational abelian factor from the individual complex character spaces.
Over the cyclotomic field, the primitive Tate module further decomposes into
two-dimensional pieces indexed by the primitive residue classes modulo $10$
\cite[Secs.~6.2--6.3]{FuselierLongRamakrishnaSwisherTu2022}; our geometric
calculation below identifies how the RM surface sits inside this primitive
Jacobian.

For the cyclic presentation $\mathcal Y_t$ in \eqref{eq:X-Y-cyclic}, let
$\xi=\zeta_{10}$ act by $y_Y\mapsto\xi y_Y$.  A convenient basis of
holomorphic character forms is
\begin{equation}\label{eq:omega-k}
 \omega_k=\frac{(1-x)^{k-1}\,dx}{y_Y^k},
 \qquad 1\le k\le9.
\end{equation}
The $\xi^{-k}$ eigenspace is one-dimensional.  Its periods satisfy the Gauss
equation
\begin{equation}\label{eq:char-Gauss}
 {}_2F_1\!\left(\frac{k}{10},1-\frac{k}{10};1;t\right).
\end{equation}
The primitive indices are
\begin{equation}\label{eq:primitive-set}
 \{1,3,7,9\},
\end{equation}
and hence the primitive factor has dimension four.

The two cyclic presentations in \eqref{eq:X-Y-cyclic} are isomorphic by
\begin{equation}\label{eq:XY-isom}
 y_X=\frac{y_Y^3}{(1-x)^2},\qquad
 y_Y=\frac{y_X^7}{x^2(1-x)^4(1-tx)^2}.
\end{equation}
Thus they are not two independent geometric families; they are two
character conventions on the same degree-ten cover.

The same curve admits the hyperelliptic model \eqref{eq:H}.  Indeed, setting
\begin{equation}\label{eq:Y-to-H}
 u=\frac{y_Y}{1-x},
 \qquad
 v=2tx-(1+u^{10})
\end{equation}
gives
\[
 u^{10}=\frac{x(1-tx)}{1-x},
\]
and hence
\[
 v^2=u^{20}+cu^{10}+1.
\]
The inverse map is
\begin{equation}\label{eq:H-to-Y}
 x=\frac{v+1+u^{10}}{2t},
 \qquad y_Y=u(1-x).
\end{equation}
Moreover
\begin{equation}\label{eq:omega-theta}
 \omega_k=-10\frac{u^{9-k}\,du}{v}.
\end{equation}

\subsection{The explicit factorization}\label{sec:factor-explicit}

The relation between the two geometries is the degree-two map
\begin{equation}\label{eq:H-D}
 p:H_t\longrightarrow D_t^- ,
 \qquad (z,w)=(u^2,uv).
\end{equation}
Its deck involution is $(u,v)\mapsto(-u,-v)$, and the cover is etale.
Composing with \eqref{eq:D-C-map} gives the degree-four map
\eqref{eq:mmap}.  The equation of $C_t^-$ follows immediately from
\[
 X^5-5X^3+5X=u^{10}+u^{-10}
\]
and \eqref{eq:H}.

The factorization is most naturally expressed by idempotent endomorphisms of
$J(D_t^-)$ in the rational isogeny category.  Here an idempotent is simply an
endomorphism $e\in\operatorname{End}^0(J(D_t^-))$ with $e^2=e$; after
multiplying by a suitable integer, its image is represented by an abelian
subvariety, well defined up to isogeny.  This language records an actual
abelian factor and not merely a chosen complex eigenspace.  Let
\begin{equation}\label{eq:idempotents}
 e_0=\frac15q^*q_*,\qquad
 e_P=1-e_0,\qquad
 e_A=e_P\frac{1+\sigma}{2},\qquad
 e_B=e_P\frac{1-\sigma}{2}.
\end{equation}
Here $e_0$ is the averaging idempotent associated with the degree-five
quotient $q:D_t^-\to E_t$, so its image is the elliptic factor.  Its
complement $e_P$ cuts out the four-dimensional connected kernel $P$ up to
isogeny.  The involution $\sigma$ then splits $P$ into its invariant and
anti-invariant abelian surfaces through the orthogonal idempotents $e_A$ and
$e_B$.  The invariant surface is the pullback of $J(C_t^-)$; hence
\begin{equation}\label{eq:P-split}
 P\sim J(C_t^-)\times B_t.
\end{equation}
On differentials one has
\begin{equation}\label{eq:pstar-diffs}
 p^*\eta_j=2\frac{u^{2j}\,du}{v}
 =-\frac{1}{5}\omega_{9-2j}.
\end{equation}
Thus $p^*P$ is precisely the primitive four-dimensional part,
and $p^*:P\to J^{\mathrm{new}}(H_t)$ is an isogeny.  The pullback from the
RM surface is especially transparent:
\begin{equation}\label{eq:m-pullback}
 m^*\left(\frac{dX}{Y}\right)=\frac{\omega_7-\omega_3}{5},
 \qquad
 m^*\left(\frac{X\,dX}{Y}\right)=\frac{\omega_9-\omega_1}{5}.
\end{equation}
This proves the first part of Theorem~\ref{thm:A}.

Over $L=\Q(\zeta_5)$ the endomorphism $r-r^{-1}$ anticommutes with
$\sigma$ and has no zero eigenvalue on $P$.  It therefore exchanges the
$e_A$- and $e_B$-parts by an isogeny, proving \eqref{eq:square}.

\smallskip
\noindent\emph{Primitive versus old quotients.}
The common order-ten generalized-Legendre cover has the additional quotient
\[
 C_t^+:\quad Y^2=(X+2)(X^5-5X^3+5X+2-4t).
\]
Its genus-two Jacobian is not interchangeable with $J(C_t^-)$.  The
minus Jacobian considered here lies in the primitive order-ten factor,
whereas the plus Jacobian lies in the order-five \emph{old} factor,
with degree-ten character indices $\{2,4,6,8\}$.  These
placements follow from the distinct character supports and the known
Frey-curve quotient correspondences
\cite[Appendix~A]{GolfieriMadriagaPacetti2025}; see also
\cite[Propositions~5.1--5.2]{ChenKoutsianasFrey2025}.
The order-twenty primitive Cohen--Wolfart factor appearing below is yet
another cyclic source; it has dimension eight, not four.

\begin{remark}[What the factorization does not say]\label{rem:factorization-scope}
The four primitive characters form a single rational Galois orbit.  Thus one
cannot obtain the RM surface by arbitrarily selecting two complex character
spaces; the extra involution is essential.  We do not determine the optimal
field of definition, degree or kernel of the isogeny between $B$ and
$J(C^-)$,
nor the full generic endomorphism algebra.  None of these refinements is
needed for the comparison developed below.
\end{remark}

\section{Integral monodromy and a modular-embedding comparison}
\label{sec:monodromy}

A representation over $K$ alone does not identify the integral homology
lattice whose reduction gives the torsion representation. Takei's
Picard--Fuchs equations and RM module supply the rank-two local system.
We fix positively oriented based loops around $0$ and $1$, choose a
basis for their unipotent monodromies, and show that the generated
$\mathcal O_K$-algebra is $M_2(\mathcal O_K)$. Class number one then
determines the integral lattice up to multiplication by a nonzero
scalar. The same matrices give the Hecke image under $\sigma_+$ and
the nondiscrete image under $\sigma_-$. They also permit comparison
with Cohen--Wolfart's Example~1 and will be reduced modulo the
coefficient primes in Section~\ref{sec:level11}.

\subsection{The integral RM lattice}\label{sec:integral-lattice}

The two equations in \eqref{eq:TTVPF} and the geometric $\OK$-module are
results of Takei \cite[Lemma~4.2.1 and Theorems~4.2.4--4.2.5]{Takei2014}.
They have nontrivial unipotent local monodromies at $0$ and $1$.  For a system
${}_2F_1(a,1-a;1;t)$ the trace at infinity is $2\cos(2\pi a)$.  If
$A,B$ are trace-$2$ unipotents with distinct fixed lines, a basis on those
fixed lines puts them into the form
\begin{equation}\label{eq:ABnormal}
 A=\begin{pmatrix}1&1\\0&1\end{pmatrix},
 \qquad
 B=\begin{pmatrix}1&0\\-\mu&1\end{pmatrix},
 \qquad
 \mu=2-\tr(AB).
\end{equation}
For the TTV pair this gives exactly $\mu=2+\alpha$ and hence the triple
\eqref{eq:exact-monodromy}.

The integrality is geometric rather than accidental.  Since
\begin{equation}\label{eq:mu-unit}
 \mu^2-3\mu+1=0,
\end{equation}
$\mu$ is a unit of $\OK$.  Consequently
\[
 M_0-I=E_{12},
 \qquad
 -\mu^{-1}(M_1-I)=E_{21}
\]
belong to the $\OK$-algebra generated by the monodromy, and therefore
\begin{equation}\label{eq:matrix-algebra}
 \OK[\rho(\pi_1(U))]=M_2(\OK).
\end{equation}
Any stable $\OK$-lattice is then $I\oplus I$ for a fractional ideal $I$.
Because $\Q(\sqrt5)$ has class number one, it is homothetic to
$\OK^2$.  This identifies the geometric lattice $\Lambda_t$, up to multiplication
by a nonzero scalar, with $\OK^2$ in the displayed matrix representation
and proves the integral assertion in Proposition~\ref{prop:B}.  It
does not identify this basis with a named symplectic basis of Takei's cycles.

\subsection{Global projective monodromy at the two real embeddings}

Let $\phi=(1+\sqrt5)/2$.  Under $\sigma_+$ one has
$\mu=\phi^2$.  Conjugating by $D=\operatorname{diag}(\phi,1)$ gives
\begin{equation}\label{eq:Hecke-generators}
 DM_0D^{-1}=T=\begin{pmatrix}1&\phi\\0&1\end{pmatrix},
 \qquad
 DM_1D^{-1}=STS^{-1},
 \qquad
 S=\begin{pmatrix}0&-1\\1&0\end{pmatrix}.
\end{equation}
This is the standard Hecke normalization \cite{LangLang2015}; the occurrence
of the Hecke geometry for this family is not new \cite{BogoLi2026}.
Writing $R=TS$, one has $R^5=-I$ and
\[
TSTS^{-1}=-R^2,\qquad (TSTS^{-1})^3=R,
\]
so the two parabolics generate the full Hecke group
$\langle S,T\rangle=\Delta(2,5,\infty)$.
The representation of the source orbifold of local signature
$(\infty,\infty,5)$ is therefore nonfaithful: its group
$\Z*C_5$ has no element of order two, whereas the image does.

Under $\sigma_-$,
\[
\mu'=\frac{3-\sqrt5}{2}\in(0,1).
\]
There is an elementary nondiscreteness proof.  Put $G_0=M_1$ in the
conjugate embedding and $G_{n+1}=G_nM_0G_n^{-1}$.  If
$G_n=\begin{psmallmatrix}a_n&b_n\\c_n&d_n\end{psmallmatrix}$, then
\[
G_{n+1}=\begin{pmatrix}1-a_nc_n&a_n^2\\-c_n^2&1+a_nc_n\end{pmatrix}.
\]
Thus $c_{n+1}=-c_n^2$ with $0<|c_0|<1$, while the $a_n$ remain bounded;
consequently $c_n\to0$ and $a_n\to1$.  The distinct matrices $G_n$
converge to $M_0$, so the conjugate image is nondiscrete.  Thus only the $\sigma_+$ representation is discrete; the shared local
projective orders do not make the two real representations discrete
triangle-group representations.

\subsection{Cohen--Wolfart: local resemblance and global difference}\label{sec:cw-comparison}

In Cohen--Wolfart's Example~1, the Schwarz functions are normalized by
$D_i(0)=P_i$, $D_i(\infty)=Q_i$, $D_i(1)=T_i$; equation~(5) gives
\begin{equation}\label{eq:CWparams}
 \left(\frac25,\frac25;\frac45\right),\qquad
 \left(\frac15,\frac15;\frac25\right),
\end{equation}
with the first component discrete
\cite[Eq.~(5) and Example~1]{CohenWolfart1990}.  For a component
${}_2F_1(s,s;2s;z)$, the change
\begin{equation}\label{eq:CWtransform}
 z=1/t,
 \qquad
 u(t)=t^{-s}y(1/t)
\end{equation}
transforms the equation to
\[
t(1-t)u''+(1-2t)u'-s(1-s)u=0.
\]
Hence the two systems become
\begin{equation}\label{eq:CWnormalized}
 \left(\frac25,\frac35;1\right),
 \qquad
 \left(\frac15,\frac45;1\right).
\end{equation}
Their order-five local projective conjugacy classes therefore match the TTV
classes crosswise.  This explains why the two constructions appear so close
at the level of local hypergeometric monodromy.

The commutator invariant \eqref{eq:kappa} separates them globally.  For the
TTV pair,
\begin{equation}\label{eq:TTV-u-polynomial}
 u^2-3u+1=0,
\end{equation}
whereas the Cohen--Wolfart pair gives
\begin{equation}\label{eq:CW-u-polynomial}
 u^2-5u+5=0.
\end{equation}
The corresponding squared-value sets $\{u^2\}$ are disjoint, so the
commutator invariant separates the two projective pairs even after Galois
conjugation permutes the two values within each orbit.  Swapping the two cusps replaces the
commutator by a conjugate of its inverse, and duality or orientation reversal
does the same; scalar twists disappear in the commutator.  Therefore the
obstruction survives all transformations relevant to the comparison and
proves Theorem~\ref{thm:C}.

\begin{remark}[The printed Cohen--Wolfart matrices]\label{rem:CWmatrix}
In Cohen--Wolfart's Example~1, pp.~104--105, the coefficient of the
second parabolic generator is printed as
\[
 \alpha=\frac{-3+\sqrt5}{2}.
\]
Taken with the adjacent generator, this value gives a nondiscrete
component, incompatible with the stated Schwarz parameters and the
assertion that this component is discrete.  An unstated normalization
or a typographical error may account for the discrepancy; we do not
propose a corrected coefficient.  Theorem~\ref{thm:C} uses their
hypergeometric parameters, Schwarz-function normalization and
equivariance relation, which are mutually consistent and do not
require this printed matrix.
\end{remark}

\subsection{Scope of the comparison}

The three constructions now occupy different positions.  Archinard provides
the general cyclic-cover construction and its primitive character decomposition.  TTV
provides the explicit genus-two RM factor whose periods realize two of these
primitive systems, together with an integral monodromy model.  Cohen--Wolfart's Example~1 provides a closely related but globally different modular embedding: the
order-five local monodromy classes overlap, but the projective representation does not.
This distinction will be important at finite level, where one must reduce the
TTV matrices themselves rather than import generators from the
Cohen--Wolfart triangle group.

\smallskip
\noindent\emph{The other Cohen--Wolfart construction.}
This negative comparison concerns the specified Example~1 on
pp.~104--105, not the Hecke cyclic-cover construction in
\cite[\S6.1, p.~106]{CohenWolfart1990}.  Both descended TTV
Hecke-period systems agree with those of this second construction.
In Section~\ref{sec:cw-hecke} we give an explicit algebraic comparison
after finite base change and determine its integral polarization.
In particular, projective
inequivalence in Theorem~\ref{thm:C} is consistent with the positive
Hecke comparison.

\section{Level-\texorpdfstring{$11$}{11} covers and quadratic Hecke geometry}
\label{sec:level11}

Once the geometric lattice $\Lambda$ is fixed, reduction at a prime
$\mathfrak p\mid11$ is a geometric operation on the family.  Analytically,
the corresponding RM torsion component is
$A_t[\mathfrak p](\C)\simeq\Lambda_t/\mathfrak p\Lambda_t$.  The
kernel of the geometric projective monodromy defines the connected
regular cover, while the degree-$12$ Borel quotient has fibre
$\mathbf P(\Lambda_t/\mathfrak p\Lambda_t)\simeq\Pone(\F_{11})$ and
records a line, equivalently a cyclic order-$11$ subgroup in
$A_t[\mathfrak p]$.  We first construct these covers directly from the integral
monodromy and distinguish the two primes of $\mathcal O_K$ above $11$.  We then compare them
with the level-$11$ congruence covers attached to the discrete Hecke group in
the sense of Clark--Voight.  The comparison is not inferred merely from a common finite
group or passport: the generator correspondence and the quadratic relation
$s=4t(1-t)$ give a connected normalized pullback of projective structures.
Because the quadratic map is ramified at the order-five point, this is not an
ordinary index-two orbifold subgroup covering.  The section concludes with
the field-of-moduli and descent statement for the Borel pairs and with the
separate classical noncongruence interpretation inside $\bar\Gamma(2)$.

\subsection{Reduction at the two primes above \texorpdfstring{$11$}{11}}

Darmon already gives the relevant projective finite realization in this
family \cite[Theorems~1.6 and~1.10 and Lemma~1.11]{Darmon2000}; we keep the ordered peripheral generators because they determine the
branch cycles of the covers below.  Modulo $11$,
\begin{equation}\label{eq:factor11poly}
 X^2+X-1=(X-3)(X-7),
\end{equation}
which gives the factorization \eqref{eq:split11}.  The reductions of
$\alpha$ are $3$ and $7$, so the reductions of $\mu=2+\alpha$ are $5$ and
$9$.  Put $A=\overline M_0$ and let $B_i,C_i$ be the reductions of
$M_1,M_\infty$ at $\mathfrak p_i$.  Then
\begin{equation}\label{eq:finite-triples}
 A=\begin{pmatrix}1&1\\0&1\end{pmatrix},
 \quad
 B_1=\begin{pmatrix}1&0\\6&1\end{pmatrix},
 \quad
 C_1=\begin{pmatrix}1&10\\5&7\end{pmatrix},
\end{equation}
while
\begin{equation}\label{eq:finite-triples2}
 B_2=\begin{pmatrix}1&0\\2&1\end{pmatrix},
 \qquad
 C_2=\begin{pmatrix}1&10\\9&3\end{pmatrix}.
\end{equation}
They satisfy
\[
AB_iC_i=I,
\]
with linear orders $(11,11,10)$ and projective orders $(11,11,5)$.

The image is large for an elementary reason.  Powers of $A$ give every
upper elementary matrix $E_{12}(x)$ and powers of $B_i$ give every lower
elementary matrix $E_{21}(y)$, $x,y\in\F_{11}$.  These generate
$\SL(\F_{11})$.  Thus
\begin{equation}\label{eq:full-images}
 G_i^{\mathrm{lin}}=\SL(\F_{11}),
 \qquad
 G_i^{\mathrm{proj}}=\PSL(\F_{11}).
\end{equation}
In particular
\[
|G_i^{\mathrm{lin}}|=1320,
\qquad
|G_i^{\mathrm{proj}}|=660.
\]

\subsection{Regular covers and their low-degree Borel quotients}

Let $X_i^{\mathrm{lin}}$ and $X_i^{\mathrm{proj}}$ be the smooth
compactifications of the regular kernel covers over the $t$-line with its ordered branch points.
Riemann--Hurwitz gives
\begin{center}
\begin{tabular}{@{}lcccc@{}}
\toprule
cover & degree & inertia at $(0,1,\infty)$ & points above & genus\\
\midrule
linear & $1320$ & $(11,11,10)$ & $(120,120,132)$ & $475$\\
projective & $660$ & $(11,11,5)$ & $(60,60,132)$ & $205$\\
\bottomrule
\end{tabular}
\end{center}
The map $X_i^{\mathrm{lin}}\to X_i^{\mathrm{proj}}$ is the quotient by
$\{\pm I\}$; it is unramified above the two parabolic branch values and
ramified with index two at the $132$ points above $t=\infty$.

The manageable object is the natural Borel quotient.  Let
$\mathcal B<\PSL(\F_{11})$ be the stabilizer of $\infty$ in the natural
action on $\Pone(\F_{11})$.  Then $|\mathcal B|=55$ and
$[G:\mathcal B]=12$.  The induced permutations for the first prime are
\begin{align}\label{eq:perm1}
\sigma_{0,1}&=(0\ 1\ 2\ 3\ 4\ 5\ 6\ 7\ 8\ 9\ 10)(\infty),\notag\\
\sigma_{1,1}&=(0)(1\ 8\ 6\ 7\ 4\ 5\ 3\ 10\ 9\ \infty\ 2),\\
\sigma_{\infty,1}&=(0\ 3\ \infty\ 9\ 1)(2)(4\ 5\ 7\ 8\ 6)(10),\notag
\end{align}
and for the second prime
\begin{align}\label{eq:perm2}
\sigma_{0,2}&=(0\ 1\ 2\ 3\ 4\ 5\ 6\ 7\ 8\ 9\ 10)(\infty),\notag\\
\sigma_{1,2}&=(0)(1\ 4\ 9\ 8\ 5\ \infty\ 6\ 3\ 2\ 7\ 10),\\
\sigma_{\infty,2}&=(0\ 7\ \infty\ 5\ 1)(2\ 10\ 4\ 6\ 8)(3)(9).\notag
\end{align}
Each triple is transitive, has product one and generates the natural faithful
copy of $\PSL(\F_{11})$ in $A_{12}$.  The common passport is exactly
\eqref{eq:passport}.  Hence
\begin{equation}\label{eq:genus3calc}
 2g(Y_i)-2=-24+10+10+8=4,
 \qquad g(Y_i)=3.
\end{equation}
These are the resulting degree-$12$ three-point covers attached to the
two TTV reductions.

\subsection{Conjugate primes but different Nielsen classes}

The nontrivial automorphism of $K/\Q$ exchanges $\mathfrak p_1$ and
$\mathfrak p_2$, but the resulting finite quotients are inequivalent as ordered branch-cycle triples.  The two
order-five elements have projective trace-squares
\begin{equation}\label{eq:trace-squares-mod11}
 \tr(C_1)^2=9,
 \qquad
 \tr(C_2)^2=5
 \quad\text{in }\F_{11},
\end{equation}
so they lie in two distinct order-five conjugacy classes.  We use the
internal labels $5A$ and $5B$ for the classes represented by $C_1$ and
$C_2$, respectively; these are not external Atlas labels.  In fact
$5B=(5A)^2$ at the level of individual conjugacy classes.  Likewise $11A$
denotes the common projective order-$11$ class represented by the parabolics.
An ordered Nielsen class here means the set of product-one generating
triples in prescribed conjugacy classes, modulo simultaneous
conjugation in $G$.  Thus the class triples are
\begin{equation}\label{eq:Nielsenclasses}
 (11A,11A,5A),
 \qquad
 (11A,11A,5B).
\end{equation}

The distinction persists even after allowing Nielsen transformations of the
free generating pair.  Indeed the canonical commutator traces are
\begin{equation}\label{eq:commutator-mod11}
 \tr[A,B_1]=5,
 \qquad
 \tr[A,B_2]=6.
\end{equation}
For determinant-one lifts this trace is independent of the signs of the
lifts, and automorphisms of $F_2$ preserve the conjugacy class of the
commutator up to inversion.  Moreover
\[
 \Aut(\PSL(\F_{11}))=\mathrm{PGL}_2(\F_{11})
\]
because the residue field is prime; see \cite[Proposition~3.2]{Kucharczyk2015}.
Thus every target-group automorphism is induced by conjugation in
$\mathrm{GL}_2(\F_{11})$, which preserves the trace of the canonical commutator
lift.  Since trace is also invariant under inversion, the values $5$ and $6$
remain distinct under all the allowed operations.  The two pairs therefore
cannot lie in the same Nielsen orbit, even after target-group automorphisms.

\begin{remark}[Finite reduction and the arithmetic Galois action]
\label{rem:GaloisCaution}
The coefficient-field automorphism exchanges the primes above $11$, but
this alone does not determine the absolute-Galois action on the Belyi
pairs. Proposition~\ref{prop:field-moduli} uses the enumeration of the
two Hecke branch-cycle classes, together with the branch-cycle Galois
action, to prove their exchange and the existence of $K$-models. The
ramification partitions alone do not supply this classification. The
two genus-three curves may still become isomorphic after their maps are
forgotten. For the regular covers, compatible arithmetic descriptions
of the deck actions and Hecke pullbacks remain to be specified.
Section~\ref{sec:cartan11} determines the constant field of the
extension of $K(t)$ fixed by the kernel of the projective arithmetic
action on $A_t[\mathfrak p]$; we call this the \emph{projective splitting
field}.
\end{remark}

\subsection{From the free-group quotient to the Hecke orbifold}

The exact monodromy gives a surjection
\begin{equation}\label{eq:F2-Delta}
 \rho:F_2=\pi_1(U)\twoheadrightarrow
 \Delta:=\Delta(2,5,\infty),
 \qquad N=\ker\rho.
\end{equation}
For each prime, let
\[
 q_i:\Delta\twoheadrightarrow G=\PSL(\F_{11}),
 \qquad
 \Delta_i=\ker q_i,
\]
and put
\[
 r_i=q_i\rho,
 \qquad
 K_i=\ker r_i.
\]
Then
\begin{equation}\label{eq:kernelidentity}
 K_i=\rho^{-1}(\Delta_i),
 \qquad
 1\longrightarrow N\longrightarrow K_i
 \longrightarrow\Delta_i\longrightarrow1,
\end{equation}
and
\begin{equation}\label{eq:quotientidentity}
 F_2/K_i\simeq\Delta/\Delta_i\simeq G.
\end{equation}
This group-theoretic identity alone does not identify the associated algebraic
covers, because the upper row comes from the ordered peripheral loops on the
three-punctured $t$-sphere, while the lower row is a congruence cover of the Hecke orbifold.
The missing geometry is supplied by the involution of the $t$-line.

Put
\[
\lambda=\frac{1+\sqrt5}{2},\qquad \mu=\lambda^2,
\]
and define
\begin{equation}\label{eq:J-involution}
 J=\begin{pmatrix}0&1\\-\mu&0\end{pmatrix}.
\end{equation}
For the matrices $T=M_0$ and $U=M_1$, direct multiplication gives
\begin{equation}\label{eq:J-equivariance}
 J^2=-\mu I,
 \qquad
 JTJ^{-1}=U,
 \qquad
 JUJ^{-1}=T.
\end{equation}
Hence $[J]$ has order two in $PGL_2(K)$.

Let $\iota(t)=1-t$ and choose the fixed base point $b=1/2$.  If
$\gamma_0$ is a positive meridian about $0$, put
$\gamma_1=\iota(\gamma_0)$ and
$\gamma_\infty=(\gamma_0\gamma_1)^{-1}$.  Then
\begin{equation}\label{eq:iota-pi1}
 \iota_*(\gamma_0)=\gamma_1,
 \qquad
 \iota_*(\gamma_1)=\gamma_0,
 \qquad
 \iota_*(\gamma_\infty)
 =\gamma_0^{-1}\gamma_\infty\gamma_0.
\end{equation}
Thus
\begin{equation}\label{eq:projective-equivariance}
 \bar\rho(\iota_*\gamma)
 =[J]\bar\rho(\gamma)[J]^{-1}
\end{equation}
for all $\gamma\in\pi_1(U,b)$.  This proves projective equivariance and
descent to the quotient orbifold.  Its invariant coordinate is
\begin{equation}\label{eq:s-map-section}
 s=4t(1-t).
\end{equation}
The points $t=0,1$ are exchanged and map unramified to the cusp $s=0$;
$t=1/2$ is fixed and maps with local degree two to the order-two point
$s=1$; and $t=\infty$ maps with local degree two to $s=\infty$.

\subsection{The quadratic hypergeometric pullback}

For $a+b=1/2$, set
\begin{equation}\label{eq:Hecke-Gauss-operator}
 \mathcal D_{a,b}
 =s(1-s)\partial_s^2+
 \left(1-\frac32s\right)\partial_s-ab.
\end{equation}
Its exponent differences at $0,1,\infty$ are
$0,1/2,b-a$.  Thus the systems
\begin{equation}\label{eq:Hecke-Gauss-systems}
 \mathcal H_1={}_2F_1\!\left(\frac3{20},\frac7{20};1;s\right),
 \qquad
 \mathcal H_2={}_2F_1\!\left(\frac1{20},\frac9{20};1;s\right)
\end{equation}
have projective local orders $(\infty,2,5)$.  The corresponding scalar
quadratic identity is classical \cite[15.8.18]{DLMF}; the occurrence of
these Hecke systems and the relation $s=4t(1-t)$ for the TTV family
also appear in \cite{BogoNikdelan2023,BogoLi2026}.  We keep
the operator calculation because it identifies the full local systems with
the chosen peripheral generators.

Let $q=t(1-t)$ and $h=1-2t$, so $s=4q$ and $h^2=1-s$.  If
$f(t)=y(s(t))$, the chain rule gives
\begin{equation}\label{eq:operator-pullback}
 \bigl(q\partial_t^2+h\partial_t-4ab\bigr)(y\circ\psi)
 =4(\mathcal D_{a,b}y)\circ\psi.
\end{equation}
The operator on the left is precisely the Gauss operator with parameters
$(2a,2b;1)$.  Therefore
\begin{equation}\label{eq:linear-pullback}
 \psi^*\mathcal H_1\simeq\mathcal L_1,
 \qquad
 \psi^*\mathcal H_2\simeq\mathcal L_2,
\end{equation}
where $\mathcal L_1,\mathcal L_2$ are the two TTV systems
\eqref{eq:TTVPF}.  The apparent singularity above $s=1$ disappears:
local solutions $A(1-s)$ and $(1-s)^{1/2}B(1-s)$ pull back to
$A(h^2)$ and $hB(h^2)$, two independent holomorphic solutions at
$t=1/2$.
\medskip
\noindent\emph{The two Cohen--Wolfart Hecke systems.}
Cohen--Wolfart's \emph{Hecke} cyclic curve in
\cite[\S6.1, p.~106]{CohenWolfart1990} has Gauss systems
\[
 (3/20,13/20;4/5),\qquad (1/20,11/20;3/5),
\]
which are distinct from the two systems in their Example~1.
For either $a=3/20$ or $a=1/20$, put
\[
 z=(1-s)^{-1},\qquad u(s)=(1-s)^{-a}y((1-s)^{-1}).
\]
Substituting into the corresponding Gauss operator yields
\[
 s(1-s)u''+\left(1-\frac32s\right)u'
    -a\left(\frac12-a\right)u=0.
\]
Thus, with algebraic scalar gauges retained, the pairs are
\[
 (3/20,13/20;4/5)\longleftrightarrow(3/20,7/20;1),
 \qquad
 (1/20,11/20;3/5)\longleftrightarrow(1/20,9/20;1).
\]
The subsequent pullback $s=4t(1-t)$ gives exactly the two TTV
systems~\eqref{eq:TTVPF}. With the peripheral convention below, these
transformations identify the two projective period maps of the
Cohen--Wolfart Hecke construction with those of the TTV family, in the
corresponding order. The fractional-power gauges alone do not identify
integral homology lattices or polarizations. The following algebraic
construction gives this comparison for a specified abelian factor.

To match the integral Hecke generators, put
\begin{equation}\label{eq:S-R-Hecke}
 S_H=\lambda^{-1}J
 =\begin{pmatrix}0&\lambda^{-1}\\-\lambda&0\end{pmatrix},
 \qquad T=M_0,
 \qquad R_H=(TS_H)^{-1}.
\end{equation}
Then
\begin{equation}\label{eq:Hecke-relations}
 S_H^2=-I,
 \qquad
 (TS_H)^5=I,
 \qquad
 S_HTS_H^{-1}=M_1,
 \qquad
 (TS_H)^2=-M_0M_1.
\end{equation}
Projectively, $\langle\bar T,\bar S_H\rangle=\Delta(2,5,\infty)$.
Choose positive meridians $a,b,c$ about
$s=0,1,\infty$ with $abc=1$ and $b^2=1$ before filling the order-five
point.  The quadratic cover has sheet character $\varepsilon(a)=0$ and
$\varepsilon(b)=\varepsilon(c)=1$; a Schreier basis for its kernel is
$a,bab^{-1}$, and the infinity meridian is $c^2$.  Therefore compatible
positive meridians satisfy
\begin{equation}\label{eq:Hecke-words}
 \gamma_0\longmapsto\bar T,
 \qquad
 \gamma_1\longmapsto\bar S_H\bar T\bar S_H^{-1}=\bar M_1,
 \qquad
 \gamma_\infty\longmapsto\bar R_H^2=\bar M_\infty.
\end{equation}
Thus the two TTV cusps are the two inverse images of the single Hecke cusp,
while the TTV order-five monodromy is the square of the Hecke order-five
monodromy.  To pass from monodromy to the projective structure, normalize
the Hecke orbifold with cusp, order-two and order-five points at
$0,1,\infty$.  The Schwarzian of the inverse uniformization and that of the
Gauss ratio for $\mathcal H_1$ have the same
double-pole coefficients $(1-\vartheta^2)/2$.
Their difference is a quadratic differential with at most simple poles at
these three points, hence vanishes because
\[
\deg K_{\Pone}^{\otimes2}(0+1+\infty)=-1.
\]
Thus the ratios agree up to a Möbius transformation, proving the
projective-structure assertion in Theorem~\ref{thm:E}, not merely equality
of local monodromy classes.

\begin{remark}[Orbifold qualification]\label{rem:orbifold-qualification}
After the order-five point at infinity is filled, the quadratic map is a
branched orbifold map, not an unramified orbifold covering.  In particular,
the degree-two sheet character of the punctured quotient does not define an
index-two subgroup of the completed Hecke triangle group.  This distinction
is essential in interpreting the level pullback below.
\end{remark}

\medskip

\subsubsection*{The integral Cohen--Wolfart--TTV comparison}
\phantomsection\label{sec:cw-hecke}
Put $E=\Q(\zeta_{20})$, $K=\Q(\sqrt5)$, and
\[
 \mathcal B_E=\operatorname{Spec}
 E[b,b^{-1},(b^{20}-1)^{-1}],\qquad
 z=b^{20},\quad T=\frac{1-b^{-10}}2,\quad c=2b^{-10}.
\]
Here $T$ denotes the TTV parameter on this finite base, so
$2-4T=c$ and $1-4T(1-T)=b^{-20}$.  Write
\[
 \mathscr X_z:\ y^{20}=x^7(1-x)^{17}(1-zx)^3,
 \qquad C_T^-:\ Y^2=X^5-5X^3+5X+c.
\]
We take the smooth projective models of these equations over
$\mathcal B_E$.  Let $A=J(C_T^-)$.
Let $d(x,y)=(x,\zeta_{20}y)$ and write
$\rho=d^*$ for the induced \emph{Picard pullback} on Jacobians.
This convention is inverse to the covariant deck action on curve homology.
The primitive factor is the abelian subscheme
\[
 A_{20}=\ker\bigl(\Phi_{20}(\rho):J(\mathscr X_z)
                  \longrightarrow J(\mathscr X_z)\bigr)^0;
\]
here $\Phi_{20}$ is the twentieth cyclotomic polynomial and the
superscript $0$ denotes the identity component.  We equip $A_{20}$
with the restriction $\lambda_{20}$ of the Jacobian principal polarization.

\begin{cwfourfold}\label{thm:cwfourfold}
There is a degree-four morphism
$f:\mathscr X_z\to C_T^-$ over $\mathcal B_E$, factoring through two
ramified double covers.  Pullback $i=f^*:A\hookrightarrow J(\mathscr X_z)$
is a closed immersion with saturated integral Betti image contained in
$A_{20}$; saturation means that the quotient of homology lattices is
torsion-free.  The four ordered deck translates define an isomorphism of
abelian schemes
\[
 \Phi:A^4\xrightarrow{\ \sim\ }A_{20},\qquad
 (a_0,a_1,a_2,a_3)\longmapsto
 \sum_{j=0}^{3}\rho^j i(a_j).
\]
Relative to the product of the canonical principal polarizations on
$A^4$, the pulled-back \emph{restricted} polarization $\Phi^*\lambda_{20}$
corresponds to the Rosati-self-adjoint matrix in $M_4(\mathcal O_K)$
\[
 H=\begin{pmatrix}
 4&0&2\varphi&0\\
 0&4&0&2\varphi\\
 2\varphi&0&4&0\\
 0&2\varphi&0&4
 \end{pmatrix},\qquad
 \varphi=\frac{1+\sqrt5}{2}.
\]
It has type $(2,2,2,2,2,2,10,10)$ and degree $40\,960\,000$.
Each single surface $i(A)$ has restricted polarization type $(4,4)$:
its transported TTV principal polarization is one quarter of this
restricted form.  Thus $\Phi$ is an isomorphism of unpolarized abelian
schemes; its pulled-back polarization is the form represented by $H$,
rather than the product polarization.
\end{cwfourfold}

\begin{proof}
\emph{The quotient maps.}
Set
\[
 F(x)=x^7(1-x)^{17}(1-zx)^3,\quad
 g(x)=\frac{1-x}{1-zx},\quad
 R(x)=\frac{1-z}{x^3(1-x)^9(1-zx)^3}.
\]
The additional automorphism
\[
 \jmath(x,y)=\bigl(g(x),R(x)y^{11}\bigr)
\]
is an involution: direct substitution gives
$F(g(x))=R(x)^{20}F(x)^{11}$ and
$R(g(x))R(x)^{11}F(x)^6=1$.  Moreover
$\jmath d\jmath=d^{11}$, so this involution does not belong to
the cyclic deck group.  In $E(b)(\mathscr X_z)$ define
\[
 \begin{aligned}
 u&=\frac{y^{10}}{x^3(1-x)^8(1-zx)^2},&
 Z&=\frac{y^{14}}{x^5(1-x)^{12}(1-zx)^2},\\
 L&=\frac{y^5}{x(1-x)^4(1-zx)},&
 W&=\frac{Z^3(u+x)}{L}.
 \end{aligned}
\]
They satisfy $u^2=xg(x)$, $Z^5u=1$, $L^2=xu$ and
\[
 \jmath(u)=u,\quad \jmath(Z)=Z,\quad
 \jmath(L)=\frac uxL,\qquad
 x^2-(1+zu^2)x+u^2=0.
\]
Consequently $W$ is invariant, and
\[
 W^2=Z(Z^{10}+2Z^5+z).
\]
The invariants $Z,W$ generate the entire fixed field: indeed
$u=Z^{-5}$, the displayed quadratic recovers $x$ up to $\jmath$,
and then
\[
 L=\frac{Z^3(u+x)}W,\qquad
 y^5=x(1-x)^4(1-zx)L,\qquad
 y^{14}=x^5(1-x)^{12}(1-zx)^2Z,
\]
recover $y=(y^5)^3/y^{14}$.  The first quotient has degree two.
Substituting $Z=b^2v$, $W=b^{11}w$ gives
\[
 D_c:\quad w^2=v(v^{10}+cv^5+1).
\]
The standard TTV involution
$\sigma(v,w)=(v^{-1},w/v^6)$ has degree-two quotient
\[
 q:D_c\longrightarrow C_T^-,\qquad
 (X,Y)=(v+v^{-1},w/v^3).
\]
The Dickson identity
$(v+v^{-1})^5-5(v+v^{-1})^3+5(v+v^{-1})=v^5+v^{-5}$
verifies the target equation.  Thus $f=q\circ\pi$, where
$\pi:\mathscr X_z\to D_c$ is the first quotient.
The genera of $\mathscr X_z,D_c,C_T^-$ are $19,5,2$;
Riemann--Hurwitz gives ramification degrees $20$ and $4$ for
$\pi$ and $q$.  Since both covers ramify, both Jacobian pullbacks
are injective (the kernel of pullback of a double cover could otherwise
define a factorizing \emph{\'etale} double cover; see
\cite[\S3]{MumfordPrym1974}).  Hence $i=f^*$ is a closed immersion,
its Betti image is saturated in the ambient Jacobian lattice, and
$f_*f^*=[4]$.

\emph{Primitive characters, four translates and the RM label.}
Write $\delta=d^2$.  On $D_c$ its action is
$(v,w)\mapsto(\zeta_{10}^{4}v,\zeta_{10}^{7}w)$.
The forms $\eta_k=v^k\,dv/w$, $0\leq k\leq4$, have
$\delta$-exponents $(7,1,5,9,3)$; the character $5$ is old.
As in \eqref{eq:m-pullback},
\[
 q^*(dX/Y)=\eta_3-\eta_1,\qquad
 q^*(X\,dX/Y)=\eta_4-\eta_0.
\]
Because $\jmath d\jmath=d^{11}$, the lifts of these two forms
to $\mathscr X_z$ have nonzero primitive $d$-character supports
\[
 M_1=(1,9,11,19),\qquad M_2=(3,7,13,17).
\]
The primitive character eigenspaces are one-dimensional in
$H^{1,0}$, by the fractional-part formula for the cyclic branch
exponents $(7,17,3,13)$.  No old character occurs, so
$i(A)\subset A_{20}$.  The symmetric correspondence
$\delta^8+\delta^{-8}$ acts on this surface as the fixed TTV RM
endomorphism $\alpha=(\sqrt5-1)/2$, as in
\cite[\S3]{TTV1991}.  This fixes the label of the two real
embeddings; the individual deck operator need not preserve $i(A)$.

In the character bases the translates $0,1,2,3$ have
Vandermonde blocks $(\zeta_{20}^{mj})_{m\in M_a,\,0\leq j\leq3}$.
Their determinants are $-4(3-\varphi)$ and $-4(2+\varphi)$,
whose product is $80\ne0$.  Since $\dim A_{20}=8$, $\Phi$ is
therefore an \emph{isogeny}.  The number $80$ is a determinant in
holomorphic character bases, \emph{not} its integral degree.

\emph{The polarization pairings between the translates.}
Let $C_n=i^\dagger\rho^ni=f_*\rho^nf^*$, where $\dagger$
denotes the adjoint for the canonical principal polarizations.
The norm formula gives $C_0=4$.  The involution $\jmath$ fixes
$i(A)$, conjugates $\rho$ to $\rho^{11}$, and $\rho^{10}=-1$
on the primitive factor.  Thus $C_n=(-1)^nC_n$.  Since homomorphism
groups of abelian varieties are torsion-free, $C_n=0$ for odd $n$.
For $C_2$, use $q^*q_*=1+\sigma^*$ to obtain
\[
 C_2=2q_*\delta^*q^*,\qquad
 q^*C_2=2\bigl(\delta^*+(\delta^*)^{-1}\bigr)q^*.
\]
On the primitive part of $J(D_c)$, write $\delta$ also for its
Picard action.  Then $\delta^5=-1$, $r=\delta^8$ and
$\delta=-r^2$.  The symmetric operator $r+r^{-1}$ restricts to
the TTV endomorphism $\alpha$ on $q^*A$, so
$\delta+\delta^{-1}=2-\alpha^2=\varphi$ there.
Injectivity of $q^*$ now gives $C_2=2\varphi$.
Rosati sends $\rho$ to $\rho^{-1}$, hence
$C_{-2}=C_2$ and the pulled-back polarization has blocks
$C_{k-j}$; this gives the displayed matrix $H$ in the order
$i,\rho i,\rho^2i,\rho^3i$.
The RM homology is free of rank two over $\mathcal O_K$ by
Proposition~\ref{prop:B}.  In the integral basis $(1,\alpha)$ of
$\mathcal O_K$, multiplication by $\varphi=1+\alpha$ has matrix
$P=\left(\begin{smallmatrix}1&1\\1&0\end{smallmatrix}\right)$.
Thus the endomorphism represented by $H$ acts on homology as two
copies of an $8$-by-$8$ integer matrix.  Reordering each copy gives
two blocks
$2\left(\begin{smallmatrix}2I_2&P\\P&2I_2\end{smallmatrix}\right)$.
The inner matrix has Smith invariants $(1,1,1,5)$.
Since the product principal polarization identifies the homology
with its dual, the full polarization homomorphism has Smith
invariants $(2^{[12]},10^{[4]})$.  Consequently
$\Phi^*\lambda_{20}$ has type $(2^6,10^2)$ and degree
$(2^6\cdot10^2)^2=40\,960\,000$.

\emph{The primitive polarization and the degree of $\Phi$.}
We compute $\deg\lambda_{20}$ independently of $\Phi$.
Wang--Xue--Yu's hypothesis~(1.1) holds with $N=20$, three
distinct finite roots, multiplicities $7,17,3$, and polynomial
degree $27$, all prime to $20$.  The change $y\mapsto y/b^3$
makes the polynomial monic over $E(b)$.  Their Theorem~1.1 and
\S3.16, equations~(3.28)--(3.30) in the cited preprint, give
\[
 c_{20}=(1-\zeta_4)^2(1-\zeta_5),\qquad
 \deg\lambda_{20}=N_{E/\Q}(c_{20})^2=6400^2.
\]
The restricted polarization is the one used here; replacing the
deck generator by its inverse does not change it.
For completeness, the following compact-homology calculation also
determines its elementary divisors.  Put
\[
 \mathcal A=\Z[t]/(t^{20}-1),\quad
 N=1+t+\cdots+t^{19},\quad
 \mathcal R=\mathcal A/(N),\quad
 p=\Phi_{20}(t),\quad g=N/p.
\]
For the four peripheral exponents $(7,17,3,13)$, the cellular
boundary sends $e_k\mapsto t^{a_k}-1$; the Fox relation column is
$e_0+t^7e_1+t^{24}e_2+t^{27}e_3$.  Compactification adds the
four cap relations $Ne_k$.  The fourth cap is redundant because
$N$ times the Fox relation equals $N\sum e_k$, and its last
coefficient is a unit.  Eliminating that relation leaves
\[
 H_1(\mathscr X_z,\Z)
 \simeq\ker\bigl[(u_7,u_{17},u_3):\mathcal R^3\to\mathcal R\bigr]
 \simeq\mathcal R^2,
 \qquad u_a=(t^a-1)/(t-1).
\]
Here $\operatorname{ann}_{\mathcal A}(t-1)=\Z N$, which justifies
passing from the compact cellular kernel to the displayed kernel.
Each $u_a$ is a unit in $\mathcal R$: if $aa'=1+20k$, then
$u_a(t)(1+t^a+\cdots+t^{(a'-1)a})=1+kN=1$ in $\mathcal R$.
Its primitive and old saturated lattices are respectively
$(g\mathcal R)^2=\ker p(t)$ and $(p\mathcal R)^2=\ker g(t)$;
these kernels are saturated because their quotients embed in a
free abelian group.  Their gluing quotient is
$\bigl(\Z[t]/(p,g)\bigr)^2$.  Evaluation at $\zeta_{20}$ gives
\[
 g(\zeta_{20})=\varepsilon\,2(1-\zeta_{20}^4),
 \qquad \varepsilon\in\Z[\zeta_{20}]^{\times},
\]
where explicit cyclotomic division gives
$\varepsilon=2t^7+2t^6-t^5-t^4+t^3+t^2-2t-2$
modulo $p$, of norm $1$.
Because $\Q(\zeta_{20})$ has degree $8$, the additive
quotient for one copy is
$(\Z/2)^8\oplus(\Z/5)^2$, and for the two copies it is
\[
 (\Z/2)^{12}\oplus(\Z/10)^4.
\]
The primitive and old rational deck spaces are orthogonal under
the unimodular Jacobian intersection form.  Their saturated
lattices are orthogonal complements.  For a primitive nondegenerate
sublattice $L$ of a unimodular lattice $\Lambda$, restriction of
the pairing gives a surjection $\Lambda\to L^*$ with kernel
$L^\perp$; hence
$\Lambda/(L\oplus L^\perp)\simeq L^*/L$.
Thus the gluing quotient is the discriminant group of the
\emph{restricted} primitive polarization.  Its alternating elementary divisors
occur twice, so $\lambda_{20}$ has type $(2^6,10^2)$ and degree
$40\,960\,000$, independently of $\Phi$.
Therefore
\[
 \deg(\Phi^*\lambda_{20})
 =(\deg\Phi)^2\deg(\lambda_{20})
\]
forces $\deg\Phi=1$.  This proves that the four surfaces sum to
the \emph{intrinsic} primitive lattice with integral index one;
the nontrivial primitive/old gluing group is not the kernel of
$\Phi$.  Finally, the branch sections $0,1,z^{-1},\infty$ stay
distinct over $\mathcal B_E$.  Their tame normalizations, the
quotient maps and all the correspondences extend over this smooth
base.  The connected primitive factor is an abelian subscheme,
obtained from the denominator-cleared cyclotomic projector with
constant character rank.  The argument applies on every geometric
fibre.  Hence $\Phi$ is proper and quasi-finite, and therefore finite;
its trivial geometric kernels and invertible differentials in
characteristic zero make it finite \emph{\'etale} of degree one.
It is an isomorphism of abelian schemes.
\end{proof}

\begin{newremark}[Relation with earlier decompositions]
Cohen--W\"ustholz already proved the fourth-power \emph{isogeny}
\cite[pp.~98--99]{CohenWustholz2002}; the primitive polarization
type follows from \cite{WangXueYu2014}, and a related integral
splitting of the genus-five dihedral component is present in
\cite[Theorem~4.8(d)]{LangeRodriguezRojas2014}.
Here the explicit degree-four quotient and its four ordered deck
translates determine the off-diagonal polarization pairings.
The integral comparison concerns these embedded subvarieties; it
does not identify a separately prescribed Pochhammer-cycle lattice or
an abelian moduli problem with specified RM action, polarization, and
level data.
\end{newremark}

\noindent\textit{Triangular congruence level.}
Clark--Voight construct congruence covers by reducing an integral order
attached to a discrete triangle group; a Borel quotient is obtained by
stabilizing a line in the resulting two-dimensional residue representation
\cite[Sec.~5 and the definition preceding (9.10)]{ClarkVoight2019}.  For the
cuspidal Hecke group used here we use the unipotent case of
Theorem~9.1 \cite[Theorem~9.1]{ClarkVoight2019}; the authors' errata explicitly
repairs this step when the corresponding triangle order is $\infty$
\cite[correction~(8)]{ClarkVoightErrata}.  The integral generators in \eqref{eq:S-R-Hecke} identify our two reductions with
these standard level-$11$ Hecke reductions, while the ordered generator
formula \eqref{eq:Hecke-words} fixes the branch marking.  Thus the curves
$E_i$ below are established triangular Borel levels; Duque--Rosero--Voight's
finite signature $(2,5,11)$ records their reduced ramification and should not
be confused with the characteristic-zero cuspidal signature
$(2,5,\infty)$ \cite{DuqueRoseroVoight2023}.

\subsection{The Hecke level at \texorpdfstring{$11$}{11}}

Let $Z_i\to\Pone_s$ be the compact regular Hecke cover with kernel
$\Delta_i=\ker q_i$, and let
\begin{equation}\label{eq:E-def}
 E_i=Z_i/\mathcal B.
\end{equation}
The notation $E_i$ means a genus-one curve; no rational origin is being
chosen.  The prime-level index is consistent with the standard Hecke-group index
formula \cite{LangLang2015}.  The regular Hecke cover has degree $660$ and
inertia $(11,2,5)$, whence
\begin{equation}\label{eq:Hecke-genus70}
 2g(Z_i)-2
 =660\left(-2+\frac{10}{11}+\frac12+\frac45\right)=138,
 \qquad g(Z_i)=70.
\end{equation}

In the natural action on $\Pone(\F_{11})$, an order-$11$ unipotent has cycle
type $11\,1$, the order-two generator is fixed-point free and has type
$2^6$, and the order-five generator has two eigenlines and type
$5^2 1^2$.  Hence
\begin{equation}\label{eq:Hecke-Borel-passport}
 \operatorname{passport}(E_i)
 =(11\,1;\ 2^6;\ 5^2 1^2),
\end{equation}
and
\begin{equation}\label{eq:E-genus-one}
 2g(E_i)-2=-24+10+6+8=0,
 \qquad g(E_i)=1.
\end{equation}

\subsection{Connected normalized quadratic base change}\label{sec:quadratic-pullback}

For each $i$, the projective TTV regular cover is the connected normalized
fibre product
\begin{equation}\label{eq:regular-fibre-product}
 \boxed{
 X_i^{\mathrm{proj}}\simeq
 \Norm\!\left(\Pone_t\times_{\Pone_s}Z_i\right)}
\end{equation}
as a $G$-cover of $\Pone_t$, geometrically over $\C$.  Indeed, on the
punctured bases the pullback monodromy is $q_i\circ\psi_*$.  The meridian
about $t=1/2$ maps to the square of the order-two Hecke generator and hence
to the identity.  Filling this puncture leaves exactly the representation already defined above
$q_i\circ\rho$ on $U$.  Since it is surjective onto $G$, the pullback is connected.  By the
monodromy--cover correspondence its compactification is a compact Riemann
surface; the resulting holomorphic map is algebraic, and the smooth model is
precisely the normalization in \eqref{eq:regular-fibre-product}.  This is a
geometric statement over $\C$, not a descent assertion.

The same argument in the transitive action on $G/\mathcal B$ gives
\begin{equation}\label{eq:borel-fibre-product}
 \boxed{
 Y_i\simeq
 \Norm\!\left(\Pone_t\times_{\Pone_s}E_i\right)}.
\end{equation}
Thus the resulting TTV Belyi curves are not merely associated with the same finite
quotients as the Hecke congruence covers; they are their normalized
quadratic base changes.

The local normalization also explains the genera.  For a point of
ramification index $e$ in a cover of the $s$-line, pullback by the quadratic
map ramified at the same base value gives $\gcd(e,2)$ normalized branches,
each with residual ramification $e/\gcd(e,2)$, while the map from the
normalized pullback to the original curve has ramification
$2/\gcd(e,2)$.  Over $s=1$ the Hecke index is $2$, so the quadratic
ramification cancels.  Over $s=\infty$ the Hecke indices are odd, so the
quadratic map remains ramified upstairs.

For the regular cover there are $660/5=132$ points over $s=\infty$, and
therefore
\begin{equation}\label{eq:genus205-from70}
 2g(X_i^{\mathrm{proj}})-2
 =2(2g(Z_i)-2)+132
 =408,
\end{equation}
recovering $g(X_i^{\mathrm{proj}})=205$.

For the Borel quotient the four points above $s=\infty$ have indices
$5,5,1,1$, all odd.  Hence $Y_i\to E_i$ is a double cover branched at
exactly four points, and
\begin{equation}\label{eq:genus3-from1}
 2g(Y_i)-2
 =2(2g(E_i)-2)+4=4,
 \qquad g(Y_i)=3.
\end{equation}
The cusp type $11\,1$ occurs independently over $t=0$ and $t=1$, the
order-two cycles disappear over $t=1/2$, and the order-five cycles are
squared without changing their lengths.  Thus the Hecke passport
\eqref{eq:Hecke-Borel-passport} becomes precisely the TTV passport
\eqref{eq:passport}.

At the level of conjugacy classes, the two Hecke order-five elements are
interchanged by squaring: with the class labels fixed above, their classes are
$5B,5A$ at $\mathfrak p_1,\mathfrak p_2$, respectively, while the TTV
infinity classes are $5A,5B$.  This agrees with the Nielsen-class distinction
proved above.

\subsection{The two genus-one Belyi pairs over \texorpdfstring{$\Q(\sqrt5)$}{Q(sqrt(5))}}\label{sec:borel-descent}

The genus-one quotients just obtained are the Borel triangular levels
predicted by the preceding congruence construction and also occur in the
enumeration of Borel-type triangular modular curves of Duque--Rosero and Voight
\cite{DuqueRoseroVoight2023}.  After cyclically ordering the branch cycles as
$(2,5,11)$, the two triples above agree with the two representatives in their
accompanying data for the entry $[2,5,11,11,11,1]$
\cite{DuqueRoseroVoightData}.  Independently, the degree-$12$ passport record
for the natural group $12T179\simeq\PSL(\F_{11})$ in the Belyi-map database
contains the same two classes
\cite{MustySchiavoneSijslingVoight2019,BelyiDBData}.  Thus the exact
permutation-class count used below is externally reproducible and is not
being inferred from the common ramification partitions alone.  The following
proposition records the arithmetic consequence of distinguishing the two
order-five classes.

\begin{proposition}[Field of moduli and descent]\label{prop:field-moduli}
Let $K=\Q(\sqrt5)$.  The two Belyi pairs
\[
 \beta_i:E_i\longrightarrow\Pone_s
\]
form one Galois orbit of size two.  Each pair has field of moduli $K$ and
admits a model over $K$.  The nontrivial automorphism of $K/\Q$ exchanges the
two isomorphism classes.  Consequently the normalized quadratic pullbacks
\[
 w^2=1-\beta_i,\qquad t=\frac{1-w}{2},
\]
and hence the genus-three covers $Y_i\to\Pone_t$, also admit models over $K$.
\end{proposition}

\begin{proof}
For the ordered passport
\[
 (11\,1;\ 2^6;\ 5^2 1^2)
\]
with monodromy group $\PSL(\F_{11})$, the two data sets cited above list
exactly two simultaneous-conjugacy classes of triples, and they are the two
classes displayed here.  They are distinguished by the two conjugacy classes
of elements of order $5$ in $\PSL(\F_{11})$.  Under the branch-cycle
action of $\Gal(\overline\Q/\Q)$, an element of order $5$ is sent to the
power prescribed by the cyclotomic character.  Exponents congruent to
$\pm1$ modulo $5$ preserve each of the two classes, whereas exponents
congruent to $\pm2$ interchange them.  The two unipotent classes of order
$11$ are interchanged by an outer automorphism of $\PSL(\F_{11})$, so they
do not introduce a further distinction.  It follows that the stabilizer of
either Belyi-pair class cuts out the quadratic subfield
$K=\Q(\sqrt5)$ of $\Q(\zeta_5)$; compare the branch-cycle discussion in
\cite[Secs.~6--7]{ClarkVoight2019}.

The automorphism group of either Belyi pair is trivial.  Indeed, an
automorphism of the cover centralizes the degree-$12$ monodromy.  It must
commute with the $11$-cycle and hence is a power of that cycle fixing its
unique fixed sheet; commuting also with the order-two permutation forces that
power to be the identity.  With trivial automorphism group, the isomorphisms between Galois
conjugates are unique and therefore satisfy Weil's cocycle condition.  This is
the standard descent criterion of Weil \cite{Weil1956}; in the Belyi setting
it is stated explicitly in \cite[Sec.~4.1]{MustySchiavoneSijslingVoight2019}.
Thus the field of moduli is a field of definition.  Finally, once
$\beta_i$ is defined over $K$, the equation $w^2=1-\beta_i$ and the map
$t=(1-w)/2$ are defined over $K$; normalization therefore gives a $K$-model
of $Y_i$.
\end{proof}

\begin{corollary}[Elliptic models over $K$]\label{cor:elliptic-K}
For each $i$, the descended genus-one curve $E_i$ has a $K$-rational point
and hence is an elliptic curve over $K$.  More precisely, the unique point of
ramification index $11$ above $s=0$ is $K$-rational and may be chosen as the
origin.
\end{corollary}

\begin{proof}
The fibre of the $K$-defined map $\beta_i$ above the $K$-rational branch
value $s=0$ is Galois stable.  It contains a unique geometric point with
ramification index $11$, so that point is fixed by
$\Gal(\overline K/K)$ and is therefore $K$-rational.
\end{proof}

\begin{corollary}[Conjugate $j$-invariants]\label{cor:j-conjugate}
Let $\tau$ be the nontrivial automorphism of $K/\Q$.  After labeling the two
pairs as in Proposition~\ref{prop:field-moduli},
\[
  j(E_2)=\tau\bigl(j(E_1)\bigr),\qquad j(E_i)\in K.
\]
Consequently
\[
 E_1\simeq_{\overline\Q}E_2
 \quad\Longleftrightarrow\quad
 j(E_1)=j(E_2)
 \quad\Longleftrightarrow\quad
 j(E_1)\in\Q.
\]
\end{corollary}

\begin{proof}
The nontrivial automorphism $\tau$ exchanges the two Belyi-pair isomorphism
classes, so it also exchanges their underlying elliptic curves.  The first
identity follows from Galois functoriality of the $j$-invariant.  The final
equivalence uses the classification of elliptic curves over an algebraically
closed field by $j$ and the fact that the fixed field of $\tau$ in $K$ is
$\Q$.
\end{proof}

\begin{remark}\label{rem:explicit-equations}
Proposition~\ref{prop:field-moduli} and
Corollaries~\ref{cor:elliptic-K}--\ref{cor:j-conjugate} are existence and
arithmetic-structure statements.  They do not supply a Weierstrass equation
for $E_i$, a rational expression for $\beta_i$, or the numerical value of
$j(E_i)$.  Explicit equations remain a separate computational problem.
\end{remark}

\subsection{The classical modular description and noncongruence}

The Hecke description is not the only uniformization of the parameter line.
Independently, the modular lambda function identifies
$\Pone_t\setminus\{0,1,\infty\}$ with
$\bar\Gamma(2)\backslash\mathfrak H$.  Under this identification the same
finite covers correspond to finite-index subgroups of $\bar\Gamma(2)$.  Their
cusp widths are obtained by multiplying the ramification indices over the
three punctures by the base cusp width $2$.  This gives generalized level
$110$ for both natural subgroups and, by the classical results of Wohlfahrt
and Larcher, proves that neither subgroup is congruence.

Fix the classical modular identification
\begin{equation}\label{eq:lambda-uniformization}
 U_t\simeq\bar\Gamma(2)\backslash\HH.
\end{equation}
The regular punctured level-$11$ cover is then uniformized by
$K_i<\bar\Gamma(2)$, while the punctured Borel quotient is uniformized by
\begin{equation}\label{eq:H-i}
 H_i=(q_i\circ\rho)^{-1}(\mathcal B).
\end{equation}
The Hecke and classical modular descriptions use different base
orbifolds.  Their relation is the normalized quadratic base change
proved above.''

\begin{corollary}[Classical noncongruence]\label{cor:noncongruence}
For each prime $\mathfrak p_i\mid11$, both the regular subgroup
$K_i<\bar\Gamma(2)$ and the Borel-preimage subgroup
$H_i<\bar\Gamma(2)$ are noncongruence subgroups of the classical modular
group.  Their generalized (Wohlfahrt) level is $110$.
\end{corollary}

\begin{proof}
Every cusp of $\bar\Gamma(2)$ has width $2$ in $PSL_2(\Z)$.  Hence a
local degree $e$ in a cover of the $t$-line contributes classical cusp width
$2e$ for the corresponding subgroup.  For $K_i$, the projective inertia
orders $(11,11,5)$ therefore give cusp widths $22,22,10$ (with the
multiplicities determined above).  For $H_i$, the degree-$12$ passport
$(11\,1;11\,1;5^2 1^2)$ gives cusp widths
\[
22,2;\qquad 22,2;\qquad 10,10,2,2.
\]
In both cases the least common multiple of all cusp widths is $110$.
By the Fricke--Wohlfahrt theorem, if either subgroup were congruence then its
classical congruence level would be $110$ \cite{Wohlfahrt1964}.  Larcher's
theorem on cusp amplitudes implies that a congruence subgroup has a cusp
whose width equals its level; a precise modern formulation is
\cite[Theorem~A]{MasonSchweizer2012}, which explicitly extends Larcher's
classical result \cite{Larcher1982}.  See also Hsu's congruence criterion
\cite[Theorem~2.3]{Hsu1996}.  Neither subgroup has a cusp of width $110$.
Hence neither is congruence.
\end{proof}

\begin{remark}[Two meanings of congruence]
Corollary~\ref{cor:noncongruence} concerns the realization inside the
classical modular group.  It does not contradict the fact that the Hecke-side
covers arise from level-$11$ congruence subgroups of the nonarithmetic triangle
group $\Delta(2,5,\infty)$.  The normalized quadratic base change relates
these two settings without identifying their notions of level.
\end{remark}

\begin{remark}[Chen's nonabelian level structures]
\label{rem:chen}
The same epimorphisms $F_2\twoheadrightarrow\PSL(\F_{11})$ can be inserted
into Chen's moduli problem for nonabelian level structures on punctured
elliptic curves \cite[Cor.~3.2.7]{Chen2018}.  Their Chen stabilizers have
index $72$ in $\SL(\Z)$, do not contain $-I$, and give genus-zero coarse
modular curves of degree $36$ over the $j$-line.
Indeed, the commutators in \eqref{eq:commutator-mod11} have projective
order $6$; these are the two order-$660$ entries with $e=6$ in
\cite[Appendix~B.2, Table~5]{Chen2018}, with the preprint numbering.
By contrast, the Borel
subgroups $H_i<\bar\Gamma(2)$ contain the central sign in their full inverse
images and give the genus-three covers $Y_i$, while the regular groups have
much larger index.  Thus these Chen components are not another moduli
description of the Borel level curves of the RM family.  We do not pursue a different
nonabelian-level interpretation here.
\end{remark}

\subsection{Arithmetic Cartan-normalizer structures at level
\texorpdfstring{$11$}{11}}\label{sec:cartan11}
We now construct the split and nonsplit Cartan-normalizer covers from the same torsion local system.
Fix $\mathfrak p=\mathfrak p_1=(11,\alpha-3)$ and work over
$K=\Q(\sqrt5)$ with $V=A_t[\mathfrak p]$ on $U_K$.
The conjugate coefficient prime $\mathfrak p_2=(11,\alpha-7)$ gives
a different residual representation.
Geometrically the residual projective image is
$G=\PSL(\F_{11})$ by \eqref{eq:full-images};
arithmetically it is $\operatorname{PGL}_2(\F_{11})$:
the self-adjoint RM decomposition makes the Weil pairing nondegenerate
on $V$, so its determinant is the cyclotomic character.
Since $K\cap\Q(\zeta_{11})=\Q$, both determinant square classes occur.
Thus the projective splitting field has constant field
$K'=K(\sqrt{-11})$, the quadratic extension cut out by the determinant
modulo squares.  The following non-Galois covers themselves are defined over $K$.

\begin{cartanlevels}\label{prop:cartanlevels}
The local system $V$ defines two geometrically connected finite
\'etale covers over $U_K$. Their $K$-models are given intrinsically by
the following torsion constructions:
$Y_{\mathrm s}$ parametrizes unordered pairs of distinct
$\F_{11}$-lines in $V$; $Y_{\mathrm{ns}}$ parametrizes
pairs $\{\ell,\ell^{(11)}\}$ of lines in
$V\otimes_{\F_{11}}\F_{121}$ not defined over $\F_{11}$, where
$\ell^{(11)}$ is obtained by the coefficient automorphism $a\mapsto a^{11}$.
Their compactifications over $\Pone_t$ have degrees $66$ and
$55$, ordered passports
\[
 \begin{array}{ll}
 Y_{\mathrm s}:&(11^6;\ 11^6;\ 1\,5^{13}),\quad g=21,\\
 Y_{\mathrm{ns}}:&(11^5;\ 11^5;\ 5^{11}),\quad g=18.
 \end{array}
\]
For a smooth parameter $t_0\in U(K)$, a $K$-point above $t_0$
is equivalent to containment of the projective arithmetic
Galois image in a conjugate of the corresponding split or
nonsplit Cartan normalizer.  This does not require either
individual residual line to be Galois-stable.
\end{cartanlevels}

\begin{proof}
The two functors are intrinsic finite Galois sets: the latter
uses Frobenius on the \emph{coefficient} field $\F_{121}$,
which commutes with arithmetic $G_K$.  Their stabilizers in
$\operatorname{PGL}_2(\F_{11})$ have orders $20$ and $24$;
intersections with $G$ have orders $10$ and $12$, so both
$G$-actions are transitive.  Using the ordered matrices
\eqref{eq:finite-triples}, the two parabolics act as
six and five $11$-cycles, respectively; the order-five element
has cycle patterns $1\,5^{13}$ and $5^{11}$.
Riemann--Hurwitz gives the stated genera; see also
\cite[Lemma~4.4]{Darmon2000}.  The torsion constructions give
the actual $K$-forms, and a fixed point in a fibre is exactly
a Galois-stable pair of the indicated type.
\end{proof}

\begin{newremark}[Specializations and quadratic quotients]
At $t=1/2$ the automorphism
$h:(X,Y)\mapsto(-X,\sqrt{-1}\,Y)$ commutes with RM, satisfies
$h^2=-1$ on the Jacobian, and has Galois conjugates $\pm h$ over $K$.
Since $u^2+1$ is irreducible over $\F_{11}$, its two eigenlines over
$\F_{121}$ form a Galois-stable unordered pair.  This gives a
$K$-point of $Y_{\mathrm{ns}}$ at the known CM specialization
\cite[\S3.2, Propositions~4--5]{TTV1991}, outside $\mathcal B_E$.

Let $v$ be a finite place of $K$ of residue characteristic outside
$\{2,5,11\}$.  If $v(t_0)>0$ or $v(1-t_0)>0$, with positive
valuation $m$ not divisible by $11$, tame inertia has projective
order $11$.  The factorization
$X^5-5X^3+5X+2=(X+2)(X^2-X-1)^2$
gives two ordinary nodes at $0$, and the degeneration at $1$ is
symmetric.  A tame generator acts as a conjugate of $M_0^m$ or $M_1^m$,
whose off-diagonal entry is nonzero modulo $\mathfrak p$.
Neither Cartan normalizer, of order $20$ or $24$, contains such
inertia; a Borel subgroup can.

The isomorphism $(X,Y)\mapsto(-X,\sqrt{-1}\,Y)$ from $C_t^-$
to $C_{1-t}^-$ has scalar Galois ambiguity and scalar square.
It therefore induces a $K$-defined involution on each projective
pair cover.  The quotients over $s=4t(1-t)$ have genera $6$
(split) and $5$ (nonsplit), obtained by the same pair actions on
$\bar T,\bar S_H,(\bar T\bar S_H)^{-1}$.
The original curves are their normalized quadratic pullbacks
$w^2=1-s$.  Away from $s=0,1,\infty$, a $K$-point downstairs
lifts if and only if $1-s\in K^{\times2}$.  At $s=1$ the
normalization must also be considered: for example, if
$1-s=u^2a(u)$ at a $K$-rational point, with $a(0)\ne0$,
the points above it are rational exactly when $a(0)$ is a square.

At the chosen coefficient prime the commutator invariant $\tr[A,B]-2$ is $3$
for the minus family and $1$ for the plus family, so the common
passports do not identify their marked covers.  The constructions
above give neither a rational-point determination nor a uniform
large-image theorem.
\end{newremark}

\section{Finite-field realization and Frobenius traces}
\label{sec:finite-field}

The preceding sections concern the Betti realization of the RM family:
first the correct abelian factor, then its integral monodromy, its explicit
Cohen--Wolfart integral comparison, and its marked level covers.
The short Cartan discussion is an arithmetic application of the same
level-$11$ local system, distinct from the finite-field Frobenius
realization considered here.  The same cyclic-cover models also carry compatible
arithmetic realizations.  The general relation between generalized-Legendre
periods, finite-field hypergeometric sums and rank-two Frobenius traces is due
to Deines--Fuselier--Long--Swisher--Tu and
Fuselier--Long--Ramakrishna--Swisher--Tu
\cite{DFLST2016,FuselierLongRamakrishnaSwisherTu2022}.  Here we match their
character conventions to the abelian factor fixed in
Section~\ref{sec:factor}.  This gives explicit Frobenius traces and the local
Euler polynomial of $A_t$.  The same calculation also supplies
an arithmetic counterpart of the quadratic Hecke transformation: normalized
finite-field hypergeometric functions agree under $s=4t(1-t)$, whereas the
unnormalized periods that equal Frobenius traces differ by the constant quadratic
character of $2$.  These statements are independent of the reduction modulo the primes
above $11$ used in Section~\ref{sec:level11}.

\subsection{Generalized-Legendre characters}

Put
\[
 L=\Q(\zeta_{10})=\Q(\zeta_5),\qquad K=L^+=\Q(\sqrt5).
\]
For Gauss parameters $(a,b;c)$, the generalized-Legendre dictionary of
\cite[(6.3)--(6.4)]{FuselierLongRamakrishnaSwisherTu2022} is
\[
 N=\operatorname{lcd}(a,b,c),\qquad
 (i,j,k)=N(1-b,1+b-c,a).
\]
Applied to the two systems in \eqref{eq:TTVPF}, it gives
\begin{equation}\label{eq:finite-field-rows}
 \left(\frac3{10},\frac7{10};1\right)
 \longleftrightarrow (10;3,7,3),
 \qquad
 \left(\frac1{10},\frac9{10};1\right)
 \longleftrightarrow (10;1,9,1).
\end{equation}
Thus the corresponding generalized-Legendre models are precisely the two
cyclic models already used in Section~\ref{sec:factor}:
\[
 \mathcal X_t:\ y_X^{10}=x^3(1-x)^7(1-tx)^3,
 \qquad
 \mathcal Y_t:\ y_Y^{10}=x(1-x)^9(1-tx).
\]
The isomorphism
\[
 y_X=\frac{y_Y^3}{(1-x)^2}
\]
sends character $1$ of $\mathcal X_t$ to character $3$ of
$\mathcal Y_t$.  With
\[
 \omega_k=(1-x)^{k-1}\frac{dx}{y_Y^k},
\]
the common-cover map of Theorem~\ref{thm:A} satisfies
\[
 m^*(dX/Y)=\frac{\omega_7-\omega_3}{5},\qquad
 m^*(X\,dX/Y)=\frac{\omega_9-\omega_1}{5}.
\]
Together with the rational projectors of Section~\ref{sec:factor}, this
identifies the two RM constituents with the inverse-character pairs
$\{3,7\}$ and $\{1,9\}$ over $L$.  The arithmetic identification below
therefore comes from the established algebraic correspondence, not merely
from agreement of Picard--Fuchs equations.

\subsection{Finite-field periods and the local Euler factor}\label{sec:frob-euler}

Every multiplicative character, including the trivial character
$\varepsilon$, is extended by zero at zero.  We use the unnormalized finite-field
period of \cite[(4.5)]{FuselierLongRamakrishnaSwisherTu2022},
\begin{equation}\label{eq:raw-P21}
 {}_2P_1[A,B;C;t;Q]
 =\sum_{x\in\F_Q}B(x)(\overline B C)(1-x)\overline A(1-tx).
\end{equation}
Let $\mathfrak q$ be a good prime of a number field containing $L$, with
residue field $\F_Q$, and let
\[
 \chi=\iota_{\mathfrak q}(1/10)
\]
be the tenth power-residue character in the convention of
\cite[Defs.~5.8--5.9]{FuselierLongRamakrishnaSwisherTu2022}.  Thus
$Q\equiv1\pmod {10}$.  We assume that the residue characteristic does not
divide $10$, that $t$ and $1-t$ are units at $\mathfrak q$, and that the
auxiliary $\ell$-adic coefficient characteristic is different from the
residue characteristic.

\begin{proposition}[Frobenius traces and the local Euler factor]
\label{prop:finite-field-euler}
Define
\begin{equation}\label{eq:TTV-finite-traces}
 a_{1,9}(t;Q)
 =-{}_2P_1[\chi,\chi^9;\varepsilon;t;Q],
 \qquad
 a_{3,7}(t;Q)
 =-{}_2P_1[\chi^3,\chi^7;\varepsilon;t;Q].
\end{equation}
These are the Frobenius traces of the two rank-two RM constituents of
$A_t=J(C_t^-)$.  They lie in $K$ and satisfy
\begin{equation}\label{eq:finite-field-conjugacy}
 a_{3,7}(t;Q)=\tau\bigl(a_{1,9}(t;Q)\bigr),
 \qquad \tau(\sqrt5)=-\sqrt5.
\end{equation}
Each constituent has determinant $Q$.  Hence, with
\[
 S=a_{1,9}+a_{3,7},\qquad N=a_{1,9}a_{3,7},
\]
the local Euler polynomial is
\begin{align}\label{eq:TTV-local-Euler}
 P_{\mathfrak q}(A_t,X)
 &=(1-a_{1,9}X+QX^2)(1-a_{3,7}X+QX^2)\notag\\
 &=1-SX+(N+2Q)X^2-QSX^3+Q^2X^4.
\end{align}
In particular,
\begin{equation}\label{eq:TTV-point-counts}
 \#C_t^-(\F_Q)=Q+1-S,
 \qquad
 \#A_t(\F_Q)=(Q+1-a_{1,9})(Q+1-a_{3,7}).
\end{equation}
\end{proposition}

\begin{proof}
For a generalized-Legendre row $(N;i,j,k)$, the inverse-deck convention of
\cite[Sec.~6.3]{FuselierLongRamakrishnaSwisherTu2022} gives the
character-$m$ Frobenius trace as
\[
 -\sum_x\chi^{-m}\!\left(x^i(1-x)^j(1-tx)^k\right)
 =-{}_2P_1[\chi^{mk},\chi^{-mi};\chi^{-m(i+j)};t;Q].
\]
For both rows in \eqref{eq:finite-field-rows}, $i+j=10$, so the lower
character is trivial.  The character relabelling and algebraic
correspondence above then identify these traces with the two constituents of $A_t$.  Opposite characters are paired by the algebraic involution,
while the coefficient automorphism $\chi\mapsto\chi^3$ exchanges the two
pairs and induces $\tau$ on $K$, giving \eqref{eq:finite-field-conjugacy}.

The real multiplication on $A_t$ is self-adjoint for the canonical
polarization.  The alternating pairing is therefore nondegenerate on each
rank-two RM constituent, and Frobenius acts on it with multiplier $Q$.
Thus each determinant is $Q$.  The trace formula of
\cite[Thm.~1.1]{FuselierLongRamakrishnaSwisherTu2022} gives the two
quadratic factors, whose product is \eqref{eq:TTV-local-Euler}.  The count
formulas follow from the zeta function of the curve and its Jacobian.
\end{proof}

The proposition uses geometric Frobenius on $H^1_{\mathrm{\acute et}}$,
equivalently arithmetic Frobenius on the covariant Tate module.  Over $L$,
Theorem~\ref{thm:A} gives
\[
 P_{\mathfrak q}(J^{\mathrm{new}}(\mathcal X_t),X)
 =P_{\mathfrak q}(A_t,X)^2;
\]
the full genus-nine cyclic Jacobian has additional old-character
contributions.  If the residue field does not contain $\mu_{10}$, one may
pass to the residue extension containing it; the character sums then compute
a power of Frobenius rather than a signed labelled trace over the smaller
field.  We make no smaller-field character formula without the corresponding
descent data.

\subsection{Arithmetic form of the quadratic Hecke transformation}\label{sec:frob-quadratic}

The normalized finite-field function of
\cite[(4.9)]{FuselierLongRamakrishnaSwisherTu2022} is
\[
 {}_2\mathbb F_1[A,B;C;t]
 =\frac{{}_2P_1[A,B;C;t]}{J(B,\overline BC)}.
\]
The distinction from the unnormalized period \eqref{eq:raw-P21} matters here.  Let
$\mathfrak q$ now have residue size $Q\equiv1\pmod {20}$, choose an
order-$20$ character $\psi$ with $\psi^2=\chi$, and put
$s=4t(1-t)$.

\begin{corollary}[Quadratic transformation of Frobenius traces]
\label{cor:finite-field-quadratic}
For $m=1,3$ and $t\notin\{0,1,1/2\}$, define
\[
 h_m(s;Q)
 =-{}_2P_1[\psi^m,\psi^{10-m};\varepsilon;s;Q].
\]
Then
\begin{equation}\label{eq:finite-field-quadratic}
 a_m(t;Q)=\eta_{\mathfrak q}\,h_m(4t(1-t);Q),
 \qquad
 \eta_{\mathfrak q}=\psi(-1).
\end{equation}
Here $a_1=a_{1,9}$ and $a_3=a_{3,7}$.  Moreover
\[
 \eta_{\mathfrak q}=(-1)^{(Q-1)/20},
\]
and this is the quadratic character of $2$ in $\F_Q^\times$.  Over
$M=\Q(\zeta_{20})$ it is the constant quadratic character attached to
$M(\sqrt2)/M$.
\end{corollary}

\begin{proof}
In \cite[(9.13)]{FuselierLongRamakrishnaSwisherTu2022}, take
$A=\psi^m$ and $B=\psi^{10-m}$.  The stated nondegeneracy hypotheses hold,
and one obtains
\[
 {}_2\mathbb F_1[\psi^m,\psi^{10-m};\varepsilon;4t(1-t)]
 ={}_2\mathbb F_1[\chi^m,\chi^{-m};\varepsilon;t].
\]
The two Jacobi denominators are
\[
 J(\chi^{-m},\chi^m)=-1,
 \qquad
 J(\psi^{10-m},\psi^{m-10})=-\psi(-1).
\]
Restoring the unnormalized periods therefore gives
\eqref{eq:finite-field-quadratic}.  Finally,
$\psi(-1)=(-1)^{(Q-1)/20}$; for $Q\equiv1\pmod {20}$ this is the
quadratic character of $2$, corresponding over $M$ to
$M(\zeta_{40})=M(\sqrt2)$.
\end{proof}

Thus the classical quadratic transformation and the projective pullback of
Section~\ref{sec:level11} have a compatible arithmetic trace relation, but
the Frobenius traces carry a constant quadratic twist that is invisible
at the level of the normalized finite-field functions.  For a fixed rational
specialization away from $0,1,1/2$, equality of these traces at all good
primes gives, by Chebotarev and Brauer--Nesbitt, the corresponding equality
of semisimplified rank-two arithmetic realizations after twisting the Hecke
realization by this quadratic character.  No integral or polarized
identification of the entire families is asserted.  This is precisely the
type of arithmetic character ambiguity that a comparison of complex
connections alone leaves open; compare \cite[Thm.~7.1]{HoffmanTu2020}.

\subsection{Examples and arithmetic interpretation}

The formulas above follow from the cited arithmetic realizations.
Direct point counts provide an independent check.  For $t=2$
and the character labelling fixed in Proposition~\ref{prop:finite-field-euler},
one obtains
\[
\begin{array}{c|c|c|c|c}
 p&a_{1,9}&a_{3,7}&P_p(A_t,X)&\eta_p\\ \hline
 41&-6-7\alpha&1+7\alpha&
 1+5X+27X^2+205X^3+1681X^4&1\\
 61&8+6\alpha&2-6\alpha&
 1-10X+102X^2-610X^3+3721X^4&-1
\end{array}
\]
with $\alpha=(\sqrt5-1)/2$ as in \eqref{eq:alpha-mu}.  Direct counts over
$\F_p$ and $\F_{p^2}$ recover the same quartic Euler polynomials; the second
row also shows that the quadratic factor in
Corollary~\ref{cor:finite-field-quadratic} is genuinely nontrivial.

For an arithmetic specialization, Proposition~\ref{prop:finite-field-euler}
computes the good-prime Frobenius traces of the two rank-two RM
constituents.  If a separate modularity theorem identifies such a constituent
with a Hilbert modular form, these are the corresponding Hecke eigenvalues in
the stated normalization; no modularity theorem is proved here.  The
Frobenius prime $\mathfrak q$, the coefficient place $\lambda\mid\ell$, and
the torsion prime $\mathfrak p_i\mid11$ of Section~\ref{sec:level11} must
remain distinct.  At a good Frobenius prime away from $11$, preservation of
a line by the arithmetic residual representation modulo $\mathfrak p_i$
forces the corresponding quadratic Frobenius polynomial to split modulo
$\mathfrak p_i$; in particular its discriminant is a square or zero.  Such
local conditions are useful tests for exceptional specializations, but
finitely many of them do not prove global reducibility, and geometric
surjectivity does not determine the arithmetic image of every fibre.

\section{Further questions}
\label{sec:outlook}

The results above settle the existence and geometric identification of the
Borel level curves, including their $K=\Q(\sqrt5)$-models, and also identify
the good-prime Frobenius traces and Euler factors of the two RM
constituents.  The remaining questions concern more explicit or more global
arithmetic refinements.

\begin{enumerate}[label=\textup{(\arabic*)},leftmargin=2.2em]
\item \textit{Explicit equations.}
Find a Weierstrass model and degree-$12$ Belyi function for each genus-one
Hecke curve $E_i$.  Once $\beta_i$ is known, the corresponding genus-three curve is
prescribed by $w^2=1-\beta_i$ and $t=(1-w)/2$; no independent genus-three
Belyi search is needed.

\item \textit{Regular covers and arithmetic markings.}
The Borel pairs and their quadratic pullbacks have $K$-models.
Section~\ref{sec:cartan11} also determines the quadratic constant
field $K'=K(\sqrt{-11})$ of the projective splitting field.
A remaining problem is to describe the regular geometric covers
with their $\PSL(\F_{11})$-deck actions over $K'$, and their descent
to $K$ allowing Galois to act on the deck group, compatibly with
the RM marking and the regular Hecke fibre-product isomorphism.

\item \textit{Specializations with reducible $\mathfrak p$-torsion.}
The Borel level curves are natural places to study specializations for which
the selected RM $11$-torsion contains a Galois-stable cyclic subgroup.
Arithmetic applications require twists, degeneration analysis and
rational-point information on these curves; geometric surjectivity of the
varying family alone gives no all-fibre big-image theorem.

\item \textit{Ordinary loci and arithmetic specializations.}
Section~\ref{sec:finite-field} gives the two good-prime Frobenius traces and
their Euler factors, but it does not classify ordinary or nonordinary fibres,
CM specializations, or exceptional residual images.  It would be useful to
compare these traces with Takei's Cartier--Manin description of the TTV
ordinary locus and with Hasse or partial-Hasse invariants.  For applications
to the Borel level curves, the local splitting conditions supplied by the
Frobenius polynomials must still be combined with global arithmetic and
rational-point information.

\item \textit{Cartan-normalizer specializations and rational points.}
The covers $Y_{\mathrm s}$ and $Y_{\mathrm{ns}}$ of
Section~\ref{sec:cartan11} have genera $21$ and $18$, with Hecke
quotients of genera $6$ and $5$.  Their full $K$-rational points,
and the square-lifting conditions from the Hecke quotients, remain
undetermined.  The known CM point at $t=1/2$ must be separated from
families with geometric endomorphism algebra exactly $K$;
fixed-level local exclusions do not imply a parameter-uniform
big-image statement.  The fourfold comparison of
Section~\ref{sec:cw-hecke} identifies the primitive
Cohen--Wolfart homology with that of the four specified
TTV translates. A comparison with a separately prescribed Pochhammer-cycle
lattice or an abelian moduli problem would require the corresponding RM
action, polarization, and level data to be specified.

\item \textit{Moduli interpretation.}
The Borel level construction used here is intrinsic to the given RM family.
A natural refinement is to ask whether it is the pullback of a standard
Hilbert-modular $\Gamma_0(\mathfrak p)$ level morphism, with the RM order,
polarization and level datum fixed precisely.  The specific
William~Y.~Chen components considered above are different curves, so such a
moduli description would require additional structure rather than a
relabelling of those components.
\end{enumerate}

A companion study treats the adjacent $(5,5,\infty)$ branch of Darmon's
programme \cite{TiwariDarmonCompanion2026}.  The literature supplies an
explicit base-change and twist relation between a TTV-form family and the
odd branch. At $r=5$, the discrete TTV constituent has image
$\Delta(2,5,\infty)$, while the discrete constituent of the companion
family has image $\Delta(5,5,\infty)$, its index-two subgroup.  Nevertheless the parameter
markings, integral representations, genera and finite covers in the two
papers remain distinct.  In particular, the ramified quadratic map of the
present paper is not the orbifold covering associated with that index-two
subgroup inclusion.


\begin{thebibliography}{BCDF25a}

\bibitem[Arc03]{Archinard2003}
N.~Archinard,
\emph{Hypergeometric abelian varieties},
Canad. J. Math. \textbf{55} (2003), 897--932.

\bibitem[BCDF25a]{BCDFI2025}
N.~Billerey, I.~Chen, L.~Dieulefait and N.~Freitas,
\emph{On Darmon's program for the generalized Fermat equation, I},
expanded version arXiv:2205.15861v5 (2025); see also the shorter published
version in J. Reine Angew. Math. (2025), DOI 10.1515/crelle-2025-0014.

\bibitem[BCDF25b]{BCDFII2025}
N.~Billerey, I.~Chen, L.~Dieulefait and N.~Freitas,
\emph{On Darmon's program for the generalized Fermat equation, II},
Math. Comp. \textbf{94} (2025), 1977--2003.

\bibitem[BL26a]{BogoLi2026}
G.~Bogo and Y.~Li,
\emph{Partial Hasse invariants for genus zero curves in Hilbert modular varieties},
arXiv:2601.16929v1 (2026).

\bibitem[BL26b]{BogoLiAtkin2026}
G.~Bogo and Y.~Li,
\emph{Atkin polynomials for families of abelian varieties with real multiplication},
arXiv:2601.16944v1 (2026).

\bibitem[BN23]{BogoNikdelan2023}
G.~Bogo and Y.~Nikdelan,
\emph{Ramanujan systems of Rankin--Cohen type and hyperbolic triangles},
Forum Math. \textbf{35} (2023), 1609--1629.

\bibitem[CK25a]{ChenKoutsianas2025}
I.~Chen and A.~Koutsianas,
\emph{Darmon's Program: A survey},
arXiv:2507.15149 (2025).

\bibitem[CK25b]{ChenKoutsianasFrey2025}
I.~Chen and A.~Koutsianas,
\emph{A modular approach to Fermat equations of signature $(p,p,5)$
using Frey hyperelliptic curves},
arXiv:2210.02316v2 (2025).


\bibitem[WCh18]{Chen2018}
W.~Y.~Chen,
\emph{Moduli interpretations for noncongruence modular curves},
Math. Ann. \textbf{371} (2018), 41--126;
preprint arXiv:1510.05687v5.

\bibitem[CV19]{ClarkVoight2019}
P.~L.~Clark and J.~Voight,
\emph{Algebraic curves uniformized by congruence subgroups of triangle groups},
Trans. Amer. Math. Soc. \textbf{371} (2019), 33--82.

\bibitem[CV26]{ClarkVoightErrata}
P.~L.~Clark and J.~Voight,
\emph{Errata and addenda: Algebraic curves uniformized by congruence subgroups of triangle groups},
online note, 6 June 2026,
\url{https://jvoight.github.io/articles/triangle-errata.pdf}.

\bibitem[CW90]{CohenWolfart1990}
P.~B.~Cohen and J.~Wolfart,
\emph{Modular embeddings for some non-arithmetic Fuchsian groups},
Acta Arith. \textbf{56} (1990), 93--110.

\bibitem[CW02]{CohenWustholz2002}
P.~B.~Cohen and G.~W\"ustholz,
\emph{Application of the Andr\'e--Oort conjecture to some questions in transcendence},
in G.~W\"ustholz (ed.),
\emph{A Panorama of Number Theory or The View from Baker's Garden},
Cambridge University Press, 2002, 89--106,
\url{https://doi.org/10.1017/CBO9780511542961.007}.


\bibitem[Dar99]{Darmon1999}
H.~Darmon,
\emph{Modularity of fibres in rigid local systems},
Ann. of Math. (2) \textbf{149} (1999), 1079--1086.

\bibitem[Dar00]{Darmon2000}
H.~Darmon,
\emph{Rigid local systems, Hilbert modular forms, and Fermat's last theorem},
Duke Math. J. \textbf{102} (2000), 413--449.

\bibitem[Dar04]{Darmon2004}
H.~Darmon,
\emph{A fourteenth lecture on Fermat's last theorem},
in \emph{Number Theory}, CRM Proc. Lecture Notes \textbf{36},
Amer. Math. Soc., Providence, RI, 2004, 103--115.

\bibitem[DFLST16]{DFLST2016}
M.~Deines, J.~Fuselier, L.~Long, H.~Swisher and F.-T.~Tu,
\emph{Generalized Legendre curves and quaternionic multiplication},
J. Number Theory \textbf{161} (2016), 175--203.

\bibitem[DV23]{DuqueRoseroVoight2023}
J.~Duque-Rosero and J.~Voight,
\emph{Triangular modular curves of small genus},
Res. Number Theory \textbf{9} (2023), Art.~3.

\bibitem[DV23d]{DuqueRoseroVoightData}
J.~Duque-Rosero and J.~Voight,
\emph{Computational data for Triangular modular curves of small genus},
\url{https://github.com/juanitaduquer/triangularModularCurves},
accessed October 2026.

\bibitem[FLRST22]{FuselierLongRamakrishnaSwisherTu2022}
J.~Fuselier, L.~Long, R.~Ramakrishna, H.~Swisher and F.-T.~Tu,
\emph{Hypergeometric functions over finite fields},
Mem. Amer. Math. Soc. \textbf{280} (2022), no.~1382.

\bibitem[GMP25]{GolfieriMadriagaPacetti2025}
I.~Golfieri Madriaga and A.~Pacetti,
\emph{Hypergeometric motives and the generalized Fermat equation},
arXiv:2412.08804v2 (2025), Appendix~A.


\bibitem[HT20]{HoffmanTu2020}
J.~W.~Hoffman and F.-T.~Tu,
\emph{Transformations of hypergeometric motives},
arXiv:2003.05031.

\bibitem[Hsu96]{Hsu1996}
T.~Hsu,
\emph{Identifying congruence subgroups of the modular group},
Proc. Amer. Math. Soc. \textbf{124} (1996), 1351--1359.

\bibitem[KraMS]{KrausManuscript}
A.~Kraus,
\emph{On the equation $x^r+y^r=z^p$},
unpublished manuscript, especially Section~V.

\bibitem[Kuc15]{Kucharczyk2015}
R.~A.~Kucharczyk,
\emph{Modular embeddings and rigidity for Fuchsian groups},
Acta Arith. \textbf{169} (2015), 77--100.

\bibitem[LL15]{LangLang2015}
C.~L.~Lang and M.~L.~Lang,
\emph{The commutator subgroup and the index formula of the Hecke group $H_5$},
J. Group Theory \textbf{18} (2015), 75--92.

\bibitem[LRR14]{LangeRodriguezRojas2014}
H.~Lange, R.~E.~Rodr\'iguez and A.~M.~Rojas,
\emph{Polarizations on abelian subvarieties of principally polarized
abelian varieties with dihedral group actions},
Math. Z. \textbf{276} (2014), 397--420.


\bibitem[Lar82]{Larcher1982}
H.~Larcher,
\emph{The cusp amplitudes of the congruence subgroups of the classical modular group},
Illinois J. Math. \textbf{26} (1982), 164--172.

\bibitem[MS12]{MasonSchweizer2012}
A.~W.~Mason and A.~Schweizer,
\emph{The cusp amplitudes and quasi-level of a congruence subgroup of $SL_2$ over any Dedekind domain},
Proc. Lond. Math. Soc. (3) \textbf{105} (2012), 311--341.

\bibitem[Mum74]{MumfordPrym1974}
D.~Mumford, \emph{Prym varieties I}, in \emph{Contributions to Analysis},
Academic Press, 1974, 325--350.


\bibitem[MSSV19]{MustySchiavoneSijslingVoight2019}
M.~Musty, S.~Schiavone, J.~Sijsling and J.~Voight,
\emph{A database of Belyi maps},
Open Book Ser. \textbf{2} (2019), 375--392.

\bibitem[BelyiDB]{BelyiDBData}
M.~Musty, S.~Schiavone, J.~Sijsling and J.~Voight,
\emph{BelyiDB computational data},
\url{https://github.com/michaelmusty/BelyiDB},
accessed October 2026.

\bibitem[DLMF]{DLMF}
NIST Digital Library of Mathematical Functions,
\emph{Hypergeometric Function, transformation formula 15.8.18},
\url{https://dlmf.nist.gov/15.8.E18}.


\bibitem[Tak14]{Takei2014}
L.~K.~Takei,
\emph{Arithmetic aspects of triangle groups},
Ph.D. thesis, McGill University, 2014.

\bibitem[TTV91]{TTV1991}
W.~Tautz, J.~Top and A.~Verberkmoes,
\emph{Explicit hyperelliptic curves with real multiplication and permutation polynomials},
Canad. J. Math. \textbf{43} (1991), 1055--1064.

\bibitem[Tiw26]{TiwariDarmonCompanion2026}
D.~Tiwari,
\emph{Hypergeometric Abelian Surfaces and Triangle-Group Levels},
manuscript, 2026.

\bibitem[WXY14]{WangXueYu2014}
H.~Wang, J.~Xue and C.-F.~Yu,
\emph{Fixed points and homology of superelliptic Jacobians},
Math. Z. \textbf{278} (2014), 169--189;
preprint arXiv:1309.0295v2.


\bibitem[Wei56]{Weil1956}
A.~Weil,
\emph{The field of definition of a variety},
Amer. J. Math. \textbf{78} (1956), 509--524.

\bibitem[Woh64]{Wohlfahrt1964}
K.~Wohlfahrt,
\emph{An extension of F. Klein's level concept},
Illinois J. Math. \textbf{8} (1964), 529--535.

\end{thebibliography}
\end{document}